\documentclass{amsart}
\usepackage{amssymb}
\usepackage{amsfonts}
\usepackage{stmaryrd}
\usepackage{mathrsfs}
\usepackage{graphicx}

\begin{document}

\newtheorem{thm}{Theorem}[section]
\newtheorem{lem}[thm]{Lemma}
\newtheorem{cor}[thm]{Corollary}
\newtheorem{pro}[thm]{Proposition}
\newtheorem{exm}[thm]{Example}

\theoremstyle{definition}
\newtheorem{defn}{Definition}[section]

\theoremstyle{remark}
\newtheorem{rmk}{Remark}[section]

\def\square{\hfill${\vcenter{\vbox{\hrule height.4pt \hbox{\vrule
width.4pt height7pt \kern7pt \vrule width.4pt} \hrule height.4pt}}}$}
\def\T{\mathcal T}

\newenvironment{pf}{{\it Proof:}\quad}{\square \vskip 12pt}

\title{The Heegaard distances cover all non-negative integers}
\author{ Ruifeng Qiu,Yanqing Zou, Qilong Guo}

%\address{School of Mathematics Science \\ Fudan University \\ Shanghai, China
%\\200433;\\
% School of Mathematics Science \\Dalian University of
%     Technology\\
%Dalian, China\\116024 }
%\email{majiming@fudan.edu.cn,qiurf@dlut.edu.cn}

%%%%%%%%%%%%%%%%%%%%%%%%%%%%%%%%%%%%%%5555

\thanks{The first author is supported  by  a grant of NSFC (No. 11171108), and the second and third authors are supported by a grant of NSFC (No. 11271058).}

%%%%%%%%%%%%%%%%%%%%%%%%%%%%%%%%%%%%%%%%%%%5

\begin{abstract}

(1) For any  integers $n\geq 1$ and $g\geq 2$, there is a closed 3-manifold $M_{g}^{n}$ which admits a distance $n$ Heegaard splitting of genus $g$ except that the pair of  $(g, n)$ is $(2, 1)$. Furthermore,   $M_{g}^{n}$ can be chosen to be hyperbolic except that the pair of $(g, n)$ is $(3, 1)$. (2) For any integers $g\geq 2$ and $n\geq 4$,  there are infinitely many non-homeomorphic closed 3-manifolds admitting distance $n$ Heegaard splittings of genus $g$.

% (2) Let $M^{*}$ be the manifold obtained by
%attaching a collection of handlebodies $\mathscr{H}$ to $\partial M$
%along a map $f$ from $\partial \mathscr{H}$ to $\partial M$. If $f$
%is a sufficiently large power of a generic pseudo-Anosov map, then
%$M^{*}$ has a Heegaard splitting of distance  $n$. The proofs rely
%essentially on Masur-Minsky's results on curve complex.
\end{abstract}

\maketitle

\vspace*{0.5cm} {\bf Keywords}: Attaching Handlebody, Heegaard distance, Subsurface Projection.\vspace*{0.5cm}


%\maketitle

\section{Introduction}

Let $S$ be a compact surface with $\chi(S)\leq -2$ not a
4-punctured sphere. Harvey \cite{h81} defined the curve complex
$\mathcal C(S)$ as follows: The vertices of $\mathcal C(S)$ are the
isotopy classes of essential simple closed curves on $S$, and $k+1$
distinct vertices $x_{0}, x_{1}, \ldots, x_{k}$ determine a
k-simplex of $\mathcal C(S)$ if and only if they are represented by
pairwise disjoint simple closed curves. For two vertices $x$ and $y$
of $\mathcal C(S)$, the distance of $x$ and $y$, denoted by $d_{\mathcal{C}(S)}(x,
y)$, is defined to be the minimal number of 1-simplexes in a
simplicial path joining x to y. In other words,
$d_{\mathcal{C}(S)}(x,
y)$ is the smallest integer $n \geq 0$ such that there is
a sequence of vertices $x_{0} = x, . . . , x_{n} =y$  such that
$x_{i-1}$ and $x_{i}$ are represented by two disjoint essential simple closed
curves on $S$  for each $1\leq i \leq n$. For two sets of vertices
in $\mathcal C(S)$, say $X$ and $Y$, $d_{\mathcal{C}(S)}(X,
Y)$ is defined to be $min\bigl\{d_{\mathcal{C}(S)}(x,
y) \
| \ x\in X, \ y\in Y\bigr\}$.  Now let $S$ be  a torus or a
once-punctured torus. In this case,  Masur and Minsky  \cite{mm99}
define $\mathcal{C}(S)$ as follows: The vertices of $\mathcal
C(S)$ are the isotopy classes of essential simple closed curves
on $S$, and $k+1$ distinct vertices $x_{0}, x_{1}, \ldots, x_{k}$
determine a k-simplex of $\mathcal C(S)$ if and only if $x_{i}$
and $x_{j}$ are represented by two simple closed curves $c_{i}$
and $c_{j}$ on $S$ such that $c_{i}$ intersects $c_{j}$ in just
one point for each $0\leq i\neq j\leq k$.
\vskip 3mm

Let $M$ be a compact orientable 3-manifold. If there is a closed surface $S$ which cuts
$M$ into two compression bodies $V$ and $W$ such that $S=\partial_{+} V=\partial_{+} W$,
then we say $M$ has a Heegaard splitting, denoted by $M=V\cup_{S} W$, where $\partial_{+}V$ (resp.$\partial _{+} W$) means the positive boundary of $V$ (resp.$W$).
%In this case, $S$ is called a Heegaard surface of $M$, and $g(S)$ is
%the genus of the splitting $M=H_{+}\cup_{S} H_{-}$. $S$ is said to be reducible if there are essential disks $D_{+}\subset H_{+}$ and $D_{-}\subset H_{-}$ such that $\partial D_{+}=\partial D_{-}$; otherwise, $S$ is said to be irreducible. $S$ is said to be weakly reducible if there are essential disks $D_{+}\subset H_{+}$ and $D_{-}\subset H_{-}$ such that $\partial D_{+}\cap\partial D_{-}=\emptyset$; otherwise, $S$ is said to be strongly irreducible. The existence of irreducible but weakly reducible Heegaard splittings have been well known. Let $g(S)\geq
%2$.
We denote by $\mathcal D(V)$(resp.$\mathcal D(W)$) the set of vertices in
$\mathcal C(S)$ such that each element of $\mathcal
D(V)$ (resp. $\mathcal D(W)$)is represented by the boundary of an essential
disk in $V$(resp. $W$).
The distance of the Heegaard splitting $V\cup_{S} W$, denoted by $d(S)$, is defined to be $d_{\mathcal{C}(S)}(\mathcal
{D}(V),
\mathcal
{D}(W))$. See \cite{h01}. \vskip 3mm

It is well known that a 3-manifold admitting a high distance Heegaard splitting has good topological and geometric properties. For example, Hartshorn\cite{ha} and
Scharlemann\cite{sch01} showed that a 3-manifold admitting a high distance Heegaard splitting contains no essential surface with small
Euler characteristic number; Scharlemann and Tomova\cite{st2} showed that a high distance Heegaard splitting is a unique minimal Heegaard splitting up to isotopy. By Geometrilization theorem and Hempel's work in \cite{h01}, a 3-manifold $M$ admitting a distance at least 3 Heegaard splitting is hyperbolic. From this view, studying Heegaard distance is an active topic on Heegaard splitting. The following is a brief survey on the existences  of high distance Heegaard splittings:\vskip 3mm

Hempel\cite{h01} showed that for any integers $g\geq 2$, and  $n\geq 2$, there is  a 3-manifold which admitting a distance at least $n$ Heegaard splitting of genus $g$. Similar results are obtained in different ways by \cite{cr} and \cite{e}. Minsky, Moriah and Schleimer \cite{mms} proved the same result for knot complements, and Li\cite{l06}constructed the non-Haken manifolds admitting high distance Heegaard splittings. In general, generic Heegaard splittings have Heegaard distances at least $n$ for any $n\geq 2$, see \cite{lm07},\cite{lm09},\cite{lm10}. By studying Dehn filling, Ma, Qiu and Zou\cite{mqz} proved that distances of genus 2 Heegaard splittings cover all non-negative integers except 1. Recently, Ido, Jang and Kobayashi\cite{ido} proved that, for any $n>1$ and $g>1$, there is a compact 3-manifold with two boundary components which admits a distance $n$ Heegaard splitting of genus $g$. And Johnson\cite{jo01} proved that there are always existing closed 3-manifold admitting a distance $n\geq 5$, genus $g$ Heegaard splitting.\vskip 3mm

The main result of this paper is the following:\vskip 2mm

{\bf Theorem 1.} For any integers $n\geq 1$ and $g\geq 2$, there is a closed 3-manifold $M_{g}^{n}$ which admits a distance $n$ Heegaard splitting of genus $g$ except that the pair of $(g, n)$ is $(2, 1)$. Furthermore,  $M_{g}^{n}$ can be chosen to be hyperbolic except that the pair of $(g, n)$ is $(3, 1)$. \vskip 3mm

{\bf Remark on Theorem 1.} (1) It is well known that there is not a distance 1 Heegaard splitting of genus 2.

(2) By the above argument, a 3-manifold admitting a distance at least 3 Heegaard splitting is hyperbolic.
Hempel \cite{h01} showed that any Heegaard splitting of a Seifert 3-manifold has distance at most 2.
Now a natural question is: For any integer $g\geq 2$, is there a closed hyperbolic 3-manifold admitting a distance 2 Heegaard splitting of genus g?
Suppose first that $g=2$. Eudave-Munoz\cite{Eu} proved that there is a hyperbolic $(1, 1)$-knot in 3-sphere, say $K$. In this case, the complement of $K$, say $M_{K}$, admits a distance 2 Heegaard splitting of genus $2$. By the main results in [1], [14] and [31], there is a slope $r$ on $\partial M_{K}$ such that the manifold obtained by doing a surgery on $M_{k}$ along $r$, say $M_{K}(r)$, is still hyperbolic. Hence $M_{K}(r)$ admits a distance 2 Heegaard splitting of genus $2$. Maybe the answer to this question has been well known when $g\geq 3$.  However we did not find  published papers related to it.

(3) If $M$ admits a distance 1 Heegaard splitting of genus 3, then $M$ contains an essential torus. Hence $M$ is not hyperbolic.

(4) The proof of Theorem 1 implies the following fact:

Let $n$ be a positive integer,  $\{F_{1},...,F_{n}\}$ be a collection of closed orientable surfaces,  and $I$ and $J=\{1,...,n\}\backslash I$ be two subsets of $\{1,...,n\}$. Then, for any integers $g\geq \max\{\sum_{i\in I} g(F_{i}), \sum_{j\in J} g(F_{j})\}$ and $m\geq 2$, there is a compact 3-manifold $M$ admitting a distance $m$ Heegaard splitting of genus $g$, say $M=V\cup_{S} W$, such that $F_{i}\subset \partial_{-} V$ for $i\in I$, and $F_{j}\subset \partial_{-}W$ for  $j\in J$. We omit the proof.\vskip 3mm

Under the arguments in Theorem 1, we have the following result: \vskip 3mm

{\bf Theorem 2.} For any integers $g\geq 2$ and $n\geq 4$,  there are infinitely many non-homeomorphic closed 3-manifolds admitting distance $n$ Heegaard splittings of genus $g$. \vskip 3mm

We organize this paper as follows:

Section 2 is devoted to introduce some results on curve complex. Then we will prove Theorem 1 when $n\neq 2$ in Section 3, Theorem 1 when $n=2$ in Section 5,  and Theorem 2 in Section 4.

\section{Preliminaries Of Curve Complex}

Let $S$ be a compact surface of genus at least 1, and $\mathcal{C}(S)$ be the curve complex of $S$.
We call a simple closed curve $c$ in $S$ is essential if $c$ bounds no disk in $S$ and is not parallel to $\partial S$. Hence each vertex of $\mathcal {C}(S)$ is represented by the isotopy class of an essential simple closed curve in $S$. For simplicity, we do not distinguish the essential simple closed curve $c$ and its isotopy class $c$ without any further notation.
The following lemma is well known, see \cite{m}, \cite{mm99}, \cite{mm00}. \vskip 2mm

{\bf Lemma 2.1.} \  $\mathcal{C}(S)$ is connected, and the diameter of $\mathcal{C}(S)$ is infinite. \vskip 2mm

We call a collection $\mathcal {G}=\{a_{0}, a_{1},...,a_{n}\}$ is a geodesic in $\mathcal {C}(S)$
if each $a_{i}\subset \mathcal {C}^{0}(S)$ and $d_{\mathcal {C}(S)}(a_{i},a_{j})=\mid i-j\mid$,
for any $0\leq i,j\leq n$. And the length of $\mathcal {G}$ is denoted by $\mathcal {L(G)}$ is defined to be $n$.
By the connection of $\mathcal {C}^{1}(S)$, there is always a shortest path in $\mathcal{C}^{1}(S)$ connecting any two vertices of $\mathcal {C}(S)$. Thus for any two distance $n$ vertices $\alpha, \ \beta$~, we call a geodesic $\mathcal {G}$ connecting $\alpha, \beta$ if $\mathcal {G}=\{a_{0}=\alpha,...,a_{n}=\beta\}$. Now for any two sub-simplicial complex $X, Y\subset \mathcal {C}(S)$, we call a geodesic $\mathcal {G}$ realizing the distance of $X$ and $Y$ if $\mathcal {G}$ connecting an element $\alpha\in X$ and an element $\beta \in Y$ such that $\mathcal {L(G)}=d_{\mathcal {C}(S)}(X,Y)$.\vskip3mm

Let $F$ be a compact surface of genus at least 1 with non-empty boundary. Similar to the definition of the curve complex $\mathcal {C}(F)$, we can define the arc and curve complex $\mathcal {AC}(F)$ as follows:

Each vertex of $\mathcal {AC}(F)$ is the isotopy
class of an essential simple closed curve or an essential properly embedded arc in $F$, and a set of vertices form a simplex of $\mathcal {AC}(F)$ if these vertices are represented by pairwise disjoint arcs or curves in $F$. For any two disjoint vertices, we place an edge between them. All the vertices and edges form 1-skeleton of $\mathcal {AC}(F)$, denoted by ${\mathcal {AC}}^{1}(F)$. And for each edge, we assign it length 1. Thus for any two vertices $\alpha$ and $\beta$ in ${\mathcal {AC}}^{1}(F)$, the distance $d_{\mathcal {AC}(F)}(\alpha, \beta)$ is defined to be the minimal length of paths in ${\mathcal {AC}}^{1}(F)$ connecting $\alpha$ and $\beta$. Similarly, we can define the geodesic in $\mathcal {AC}(F)$.\vskip 3mm

When $F$ is a subsurface of $S$, we call $F$ is essential in $S$ if the induced map of the inclusion from
$\pi_{1}(F)$ to $\pi_{1}(S)$ is injective. Furthermore, we call $F$ is a
proper essential subsurface of $S$ if $F$ is essential in $S$ and at least one boundary component of $F$ is essential in $S$. For more details, see \cite{mm00}.  \vskip 3mm

So if $F$ is an essential subsurface of $S$, there is some connection between the $\mathcal {AC}(F)$ and $\mathcal {C}(S)$. For any $\alpha\in \mathcal {C}^{0}(S)$, there is a representative essential simple closed curve $\alpha_{geo}$ such that the intersection number $i(\alpha_{geo}, \partial F)$ is minimal. Hence each component of $\alpha_{geo}\cap F$ is essential in $F$ or $S-F$. Now for $\alpha\in \mathcal {C}(S)$, let $\mathcal {\kappa}_{F}(\alpha)$ be isotopy classes of the essential components of $\alpha_{geo}\cap F$.
It is well defined \textbf{$\triangleright \triangleright$}since for any two isotopy class $\alpha_{1}$ and $\alpha_{2}$ of $\alpha$ which both intersect $\partial F$ minimally, either

(1) $\alpha_{1}\cap \alpha_{2}=\emptyset$. Then they bounds an annulus $A$ in $S$. Hence either

(1.1) $\alpha_{1}\cap \partial F=\emptyset$.
If $\alpha_{1}$ is essential in $F$, then $A\cap F=\emptyset$. Hence $A\subset F$. And $\alpha_{1}\cap F$ is isotopic to $ \alpha_{2}\cap F$ in $F$.
If both $\alpha_{1}\cap F$ and $\alpha_{2}\cap F$ are inessential, then the essential components of
$\alpha_{1}\cap F$ and $\alpha_{2}\cap F$ are $\emptyset$. Or,

(1.2) $\alpha_{1}\cap \partial F\neq \emptyset$.

By minimality of intersection number, $A\cap \partial F$ consists of squares. It is not hard to see that
each component $\alpha_{1}\cap F$ (resp. $\alpha_{2}\cap F$) is essential. And each component $\alpha_{1}\cap F$ is
isotopic to one component of $\alpha_{2}\cap F$. The reverse is true. Or,

(2) $\alpha_{1}\cap \alpha_{2}\neq \emptyset$.
Since the intersection number $i(\alpha_{1},\alpha_{2})=0$, by Bigon Criterion (proposition 1.7\cite{FM}), there is always an innermost Bigon $B$ bounded by $\alpha_{1}\cup \alpha_{2}$
in $S$. Since there is no proper bigon in $B$ bounded by $\partial B $ and $\partial F$ (the minimality of $\alpha_{1}\cap \partial F$ and $\alpha_{2}\cap \partial F$), we can isotopy $\alpha_{1}$ and $\alpha_{2}$ such that $\alpha_{1}^{'}\cup \alpha$ (resp. $\alpha_{2}^{'}\cup \alpha_{2}$) bounds an annulus in $S$ and $\alpha_{1}^{'}\cap \partial F$ (resp. $\partial F\cap \alpha_{2}^{'}$) is minimal. And the Bigon $B$ is vanishing. Following (1), we get that any essential component of $\alpha_{1}\cap F$(resp.
$\alpha_{2}\cap F$) is isotopic to an essential component of $\alpha_{1}^{'}\cap F$ (resp.$\alpha_{2}^{'}\cap F$). And the reverse is true. We can do it again and again until there is no such Bigon. Then it returns to (1).\textbf{$\triangleleft \triangleleft $}

For any $\gamma\in \mathcal {C}(F)$, $\gamma'\in \mathcal {\sigma}_{F}(\beta)$ if and only if $\gamma'$ is the essential boundary component of a closed regular neighborhood of $\gamma\cup \partial {F}$. Specially, let $\mathcal {\sigma}_{F}(\emptyset)=\emptyset$. Now let $\pi_{F}= {\mathcal {\sigma}_{F}}\circ {\mathcal {\kappa}_{F}}$. Then the map $\pi_{F}$ links the $\mathcal {AC}(F)$ and $\mathcal {C}(S)$, which is the defined subsurface projection map in \cite{mm00}.\vskip 2mm

We say $\alpha\in \mathcal {C}^{0}(S)$ cuts $F$ if $\pi_{F}(\alpha)\neq \emptyset$. If $\alpha$, $\beta\in \mathcal {C}^{0}(S)$ both cut $F$, we write $d_{\mathcal {C}(F)}(\alpha, \beta)= diam_{\mathcal {C}(F)}(\pi_{F}(\alpha),\pi_{F}(\beta))$. And if $d_{\mathcal {C}(S)}(\alpha, \beta)=1$, then $d_{\mathcal {AC}(F)}(\alpha, \beta)\leq 1$ and $d_{\mathcal {C}(F)}(\alpha, \beta)\leq 2$, observed by H.Masur and Y.N.Minsky at first. What if the two vertices $\alpha$ and $\beta$ has distance $k$ in $\mathcal {C}(S)$? \vskip 2mm

The following is immediately followed from the above observation.

{\bf Lemma 2.2.} Let  $F$ and $S$ be as above, $\mathcal {G}=\{\alpha_{0},\ldots, \alpha_{k}\}$ be a geodesic of $\mathcal {C}(S)$ such that $\alpha_{j}$ cuts $F$ for each $0\leq i\leq k$. Then $d_{\mathcal {C}(F)}(\alpha_{0}, \alpha_{k})\leq 2k$. \vskip 2mm

However, when $k$ is quite large, the Lemma 2.2 can not provide more information. In general,  Masur-Minsky  \cite{mm00} proved the following result called Bounded Geodesic Image Theorem. \vskip 2mm

{\bf Lemma 2.3.} Let $F$ be an essential sub-surface of $S$,
 and $\gamma$ be a geodesic segment in
$\mathcal {C}(S)$, such that $\pi_{F}(v)\neq \emptyset$ for every
vertex $v$ of $\gamma$. Then there is a constant $\mathcal{M}$
depending only on $S$ so that $diam_{\mathcal {C}({F})}
(\pi_{F}(\gamma))\leq \mathcal{M}$. \vskip 3mm

When $S$ is closed with $g(S)\geq 2$, there is always a compact 3-manifold $M$ with $S$ as its compressible boundary. Let $\mathcal{D}(M, S)$, called disk set for $S$, be the subset of vertices of $\mathcal C(S)$, where each element bounds a disk in $M$. Now an essential simple closed curve on $S$, say $c$,  is said to be disk-busting if $S-c$ is incompressible in $M$. Since any two essential disks intersect in a typical way, it provides more information to study the subsurface projection of disk complex. The following Disk Image Theorem is proved by T.Li \cite{l}, H.Masur and S.Schleimer \cite{ms} independently. \vskip 3mm

{\bf Lemma 2.4.}  Let $M$ be a compact orientable  and irreducible 3-manifold. S is a boundary component of $M$.
Suppose $\partial M -S$ is incompressible. Let $\mathcal {D}$
be the disk complex of $S$, $F\subset S$ be an essential subsurface. Assume each component of $\partial F$
is disk-busting. Then either

(1) M is an I-bundle over some compact surface, $F$ is a horizontal boundary of the I-bundle
and the vertical boundary of this I-bundle is a single annulus. Or,

(2) The image of this complex, $\mathcal {\kappa}_{F}(\mathcal {D})$, lies in a ball of radius 3 in
$\mathcal {AC}(F)$. In particular, $\mathcal {\kappa}_{F}(\mathcal {D})$ has diameter 6 in $\mathcal {AC}(F)$.
Moreover, $\pi_{F}(\mathcal {D})$ has diameter at most 12 in $\mathcal {C}(F)$.\vskip 3mm

{\bf Note.} For any I-bundle ${J}$  over a bounded compact surface $P$, $\partial J= \partial_{v}{J}\cup \partial_{h}{J}$, where
the vertical boundary $\partial_{v}{J}$ is the I-bundle related to $\partial P$, and the horizontal boundary $\partial_{h}{J}$ is the portion of $\partial J$ transverse to the I-fibers.\vskip 3mm

On the other side, J.Hempel \cite{h01} defined a full simplex $X$ on $S$ to be a dimension $3g(S)-4$ simplex in $\mathcal {C}(S)$. Hence, after attaching 2-handles and 3-handles along the vertices of $X$ in the same side of  $S$ from the same side, we can get a handlebody, denoted by $H_{X}$.\vskip 3mm

{\bf Lemma 2.5 \cite{h01}.} Let $S$ be a  closed, orientable surface of genus at least 2. For any positive number $d$, any full simplex $X$ of $\mathcal {C}(S)$, there is another  full simplex $Y$ of $\mathcal {C}(S)$ such that $d_{\mathcal {C}(S)}(\mathcal {D}(H_{X}), \mathcal {D}(H_{Y}))\geq d$.\vskip 3mm

Through subsurface projection, the Bounded Geodesic Image theorem links the geodesic in curve complex and a proper subsurface. The example 1.5 \cite{mm00} shows that there is a geometry
rigidity in curve complex. With the property of infinity of diameter of curve complex, we can construct any long geodesic in curve complex. Furthermore, we also require that the constructed geodesic satisfying some condition, such as the first and last vertices are represented
by separating essential simple closed curves.

We organize our results as the following lemma which is a more stronger version of Lemma 4.1 in \cite{mqz}. \vskip 2mm

{\bf Lemma 2.6.} \ Let $g, n, m, s, t$ be be five integers such that $g, m, n\geq 2$, and  $1\leq t,s\leq g-1$. Let $S_{g}$ be a closed surface of genus $g$. Then there are two essential separating curves $\alpha$ and $\beta$ in $S_{g}$ such that $d_{\mathcal{C}(S_{g})}(\alpha, \beta)=n$, one component of $S_{g}-\alpha$ has genus $t$ while one component of $S_{g}-\beta$ has genus $s$. Furthermore,  there is a geodesic $\mathcal {G}=\{a_{0}=\alpha, a_{1}, ..., a_{n-1}, a_{n}=\beta\}$ in $\mathcal {C}(S_{g})$ such that

(1)  $a_{i}$ is non-separating in $S_{g}$ for  $1\leq i\leq n-1$, and

(2) $m\mathcal{M}+2\leq d_{\mathcal {C}(S^{a_{i}})}(a_{i-1}, a_{i+1})= m\mathcal{M}+6$,
where $S^{a_{i}}$ is the surface $S-N(a_{i})$ for $ 1\leq i \leq n-1$. \vskip 3mm

{\bf Proof.} Let $\alpha$ be an essential separating curve in $S$ such that one component of $S_{g}-\alpha$, say $S_{1}$, has genus $t$.

%When n=1, for $g\geq 3$, we can choose a separating slope $\beta$ in $S_{2}$ such that
%one component of $S-N(\beta)$ has genus $s$. It is not hard to see that $\beta$ is separating and disjoint
%from $\alpha$.

Suppose first that $n=2$. Let $b$ be a non-separating curve in $S_{g}$ which is disjoint from $\alpha$. Let$S^{b}$  be the surface $S_{g}-N(b)$, where $N(b)$ is a open regular neighborhood of $b$ in $S_{g}$. Then $S^{b}$ is a genus $g-1$ surface with two boundary components. Furthermore, $\alpha$ is an essential separating simple closed curve in $S^{b}$.

By Lemma 2.1, ${\mathcal {C}^{1}(S^{b})}$ is connected and its diameter is infinite. Hence there is an essential simple closed curve $c$ in $S^{b}$ with ${d_{\mathcal C(S^{b})} (\alpha, c)}=m\mathcal{M}+4$. Note that $g-1\geq 1$. If $c$ is separating in $S^{b}$,  there is a non-separating essential simple closed curve
$c^{*}$ in $S^{b}$ such that $c\cap c^{*}=\emptyset$. Hence ${d_{\mathcal C(S^{b})} (c, c^{*})}=1$, and $m\mathcal {M}+3\leq d_{\mathcal C_{S^{b}}}(\alpha, c^{*})\leq m\mathcal{M}+5$. It means there is a non-separating slope $c$ in $S^{b}$ such that $m\mathcal {M}+3\leq d_{\mathcal C_{S^{b}}}(\alpha, c)\leq m\mathcal{M}+5$.

Let $l$ be a non-separating simple closed curve in $S^{b}$ such that $l$ intersects $c$ in one point, and $e$ be the boundary of the regular neighborhood of $c\cup l$. Then $e$ bounds a once-punctured torus $T$ containing $l$ and $c$. Since $s\leq g-1$, there is an essential separating simple closed curve $\beta$ such that $\beta$ bounds a once-punctured surface of genus $s$ containing $T$ as a sub-surface. See Figure 1.

\begin{center}
\includegraphics[height=3cm, width=8cm]{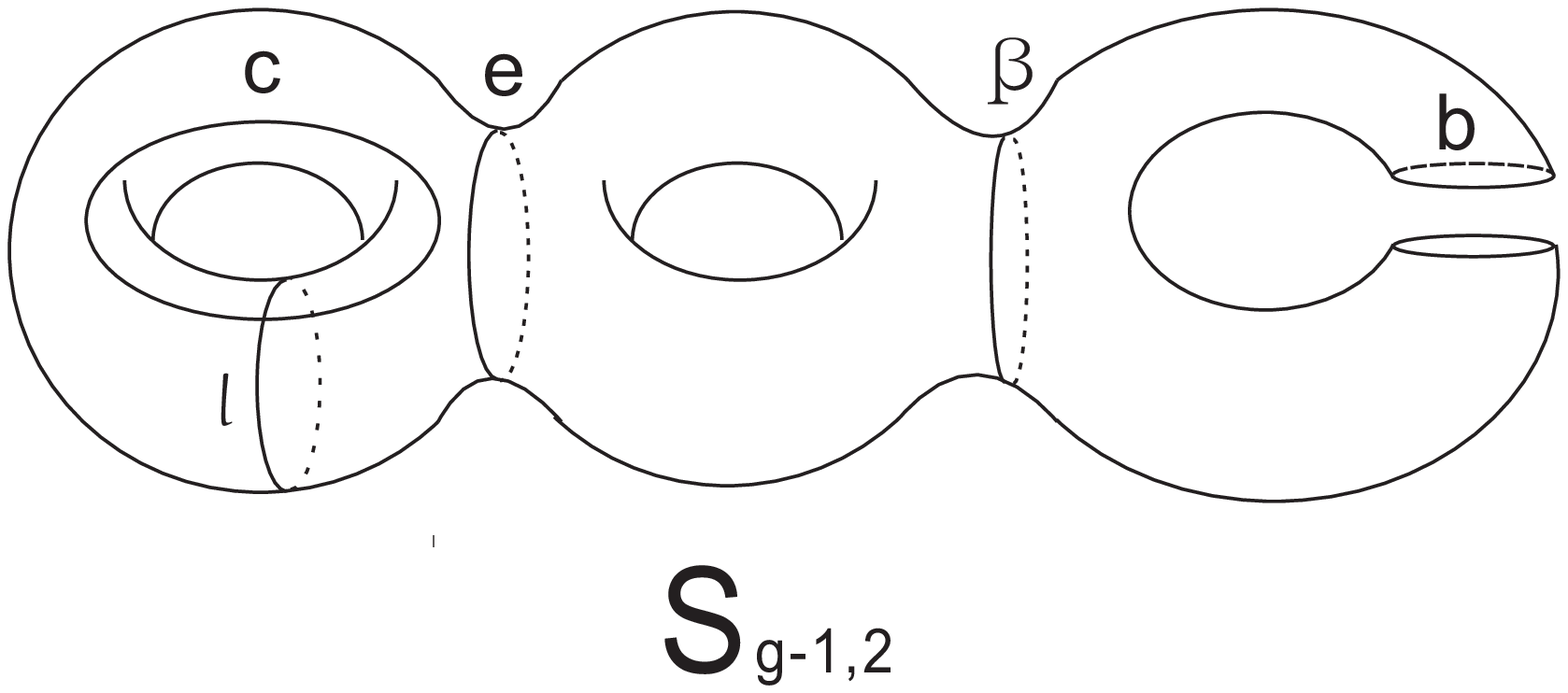}
\begin{center}
Figure 1
\end{center}
\end{center}\vskip 2mm

It is easy to see that $\beta$ is also separating in $S_{g}$. Now we prove that $d_{\mathcal{C}(S_{g})}(\alpha, \beta)=2$, and ${d_{\mathcal{C}(S_{g})}(\alpha, c)=2}$.

Since $\alpha\cap b=\emptyset$, $\beta\cap b=\emptyset$ and $c\cap b=\emptyset$, ${d_{\mathcal{C}(S_{g})}(\alpha,\beta)\leq 2}$, and ${d_{\mathcal{C}(S_{g})}(\alpha, c)\leq 2}$.
Since $c\cap \beta=\emptyset$, by the assumption on ${d_{\mathcal C(S^{b})} (\alpha, c)}$, $m\mathcal {M}+2\leq d_{\mathcal {C}(S^{b})}(\alpha, \beta)\leq m\mathcal{M}+6$.
Hence ${d_{\mathcal {C}(S_{g})}(\beta,\alpha)=2}$. For if ${d_{\mathcal{C}(S_{g})}(\alpha,\beta)\leq 1}$, then, by Lemma 2.3, $d_{\mathcal {C}(S^{b})}(\alpha, \beta)\leq \mathcal{M}$, a contradiction. Similarly, $d_{\mathcal{C}(S_{g})}(\alpha, c)=2$. And $ \mathcal {G}=\{a_{0}=\alpha, a_{1}=b, a_{2}=\beta\}$ and $\mathcal {G}^{*}=\{a_{0}=\alpha, a_{1}=b, a_{2}=c\}$ are two geodesics of $\mathcal{C}(S_{g})$.  Furthermore, $\mathcal {G}$ satisfies the conclusion of Lemma 2.6.\vskip3mm

Now we prove Lemma 2.6 by induction on $n$. \vskip 3mm

{\bf Assumption.} \ Let $k\geq 2$. Suppose that there are two essential separating simple closed curves $\alpha$ and $\beta$, and a non-separating simple closed curve $c$ in $S_{g}$ such that  ${d_{\mathcal {C}(S_{g})}(\alpha,\beta)=k}$, ${d_{\mathcal {C}(S_{g})}(\alpha,c)=k}$, and one component of $S_{g}-\alpha$ has genus $t$ while one component of $S_{g}-\beta$ has genus $s$. Furthermore, there is a geodesic $\mathcal {G}=\{\alpha=a_{0}, a_{1}, ..., a_{k-1}, a_{k}=\beta\}$ satisfying Lemma 2.6(1) and (2), and a geodesic $\mathcal{G}^{*}=\{\alpha, a_{1},...,a_{k-1}, c\}$ is also a geodesic connecting $\alpha$ to $c$ satisfying $m\mathcal {M}+3\leq d_{\mathcal {C}(S^{a_{i}})}(a_{i-1},a_{i+1})\leq m\mathcal {M}+5$, for any $1\leq i\leq k-2$ and $m\mathcal {M}+3\leq d_{\mathcal {C}(S^{a_{k-1}})}(a_{k-2},c)\leq m\mathcal {M}+5$,

%such that $t\mathcal {M}+3\leq d_{\mathcal {C}(S_{a_{1}^{*}})}(a_{i-1}^{*},a_{i+1}^{*})\leq t\mathcal {M}+5$ where $S_{a_{i}^{*}}=S-N(a_{i}^{*})$.\vskip 3mm

%When $n=K+1$.
Let $S^{c}$ be the surface $S_{g}-N(c)$, where $N(c)$ is a open regular neighborhood of $c$ in $S_{g}$. Since $c$ is non-separating in $S_{g}$, $S^{c}$ is a genus $g-1$ surface with two boundary components. Since $\mathcal{G}^{*}=\{\alpha, a_{1},...,a_{k-1}, c\}$ is also a geodesic connecting $\alpha$ to $c$, $a_{k-1}$ is an essential non-separating simple closed curve in $S^{c}$.
By the above argument, there are  a non-separating curve $h$ and a separating curve $e$ in $S^{c}$ such that $e$ bounds a once-punctured torus $T^{*}$ containing $h$, $m\mathcal {M}+3\leq d_{C(S^{c})} (h, a_{k-1})\leq m\mathcal {M}+5$, and $m\mathcal {M}+2\leq {d_{C(S^{c})}(e,a_{k-1})}\leq m\mathcal {M}+6$. And there is aslo an essential separating simple closed curve $\gamma$ which bounds a genus $s$ sub-surface of $S^{c}$ containing $T^{*}$ as a sub-surface. Not hard to see $\gamma$ is also separating in $S_{g}$. Since $h$ is disjoint from $\gamma$,  $m\mathcal {M}+2 \leq d_{C(S^{c})}(\gamma, a_{k-1})\leq m\mathcal {M}+6$. \vskip3mm

Now we prove that  $d_{\mathcal{C}(S_{g})}(\alpha,h)=k+1$, $d_{\mathcal{C}(S_{g})}(\alpha,\gamma)=k+1$. \vskip3mm

Suppose,otherwise, that $d_{\mathcal{C}(S_{g})}(\alpha,h)=x\leq k$. Then there exists a geodesic line ${\mathcal {G}_{1}}=\bigl\{\alpha=b_{0},\ldots, b_{x}=h\bigr\}$. Note that each of $\alpha$ and $h$ is not isotopic to $c$. with the length is less then or equal to $K$. Since ${d_{\mathcal {C}(S_{g})}(\alpha,c)=k}$, $b_{j}$ is not isotopic to $c$ for $1\leq j\leq x-1$. That means $b_{j}$ cuts $S^{c}$ for each $0\leq j\leq x$. By Lemma 2.3, ${d_{\mathcal{S}^{c}}(\alpha, h)\leq \mathcal {M}}$. On the other side, since ${d_{\mathcal {C}(S_{g})}(\alpha,c)=k}$, $a_{j}$ is not isotopic to $c$ for $0\leq j\leq k-1$. By using Lemma 2.3 again, ${d_{\mathcal{S}^{c}}(\alpha, a_{k-1})\leq \mathcal {M}}$. Then ${d_{\mathcal{C}({S}^{c})}(a_{k-1}, h)\leq 2\mathcal {M}}$. It contradicts the choice of $h$.
\vskip 3mm

Now $\mathcal {G}^{'}=\{a_{0}=\alpha, a_{1},...,a_{k-1}, c, \gamma\}$ and $\mathcal {G}^{''}=\{a_{0}=\alpha, a_{1},...a_{k-1}, c,  h\}$ are two geodesics satisfying the Assumption.
Hence Lemma 2.6 holds. END.\vskip 3mm

\section{Proof of Theorem 1 (1)}

In this section, we will prove  the following proposition: \vskip 2mm

{\bf Proposition 3.1.} \ For any  positive integers $n\neq 2$ and $g\geq 2$, there is a closed 3-manifold which admits a distance $n$ Heegaard splitting of genus $g$ except that the pair of $(g, n)$ is $(2, 1)$. Furthermore, $M_{g}^{n}$ can be chosen to be hyperbolic except that the pair of $(g, n)$ is $(3, 1)$. \vskip 2mm

{\bf Proof.} \ We first suppose that $n\geq 3$. \vskip 3mm

Let $S$ be a closed surface of genus $g$. By Lemma 2.6, there are two separating essential simple closed curves $\alpha$ and $\beta$ such that $d_{\mathcal {C(S)}}(\alpha, \beta)=n$ for $n\geq 3$. Let $V$ be the compression body obtained by attaching a 2-handle to $S\times [0,1]$  along $\alpha\times \bigl\{1\bigr\}$, and $W$ be the compression body obtained by attaching a 2-handle to $S\times [-1,0]$ along $\beta\times \bigl\{-1\bigr\}$. Then $V\cup_{S} W$ is a Heegaard splitting where $S$ is the surface $S\times\bigl\{0\bigr\}$, see Figure 2. Since $V$ contains only one essential disk $B$ with $\partial B=\alpha$ up to isotopy, and $W$ contains only one disk $D$ with $\partial D=\beta$ up to isotopy, $d_{\mathcal {C}(S)}(V,W)=n$.

\begin{center}
\includegraphics[height=6cm, width=10cm]{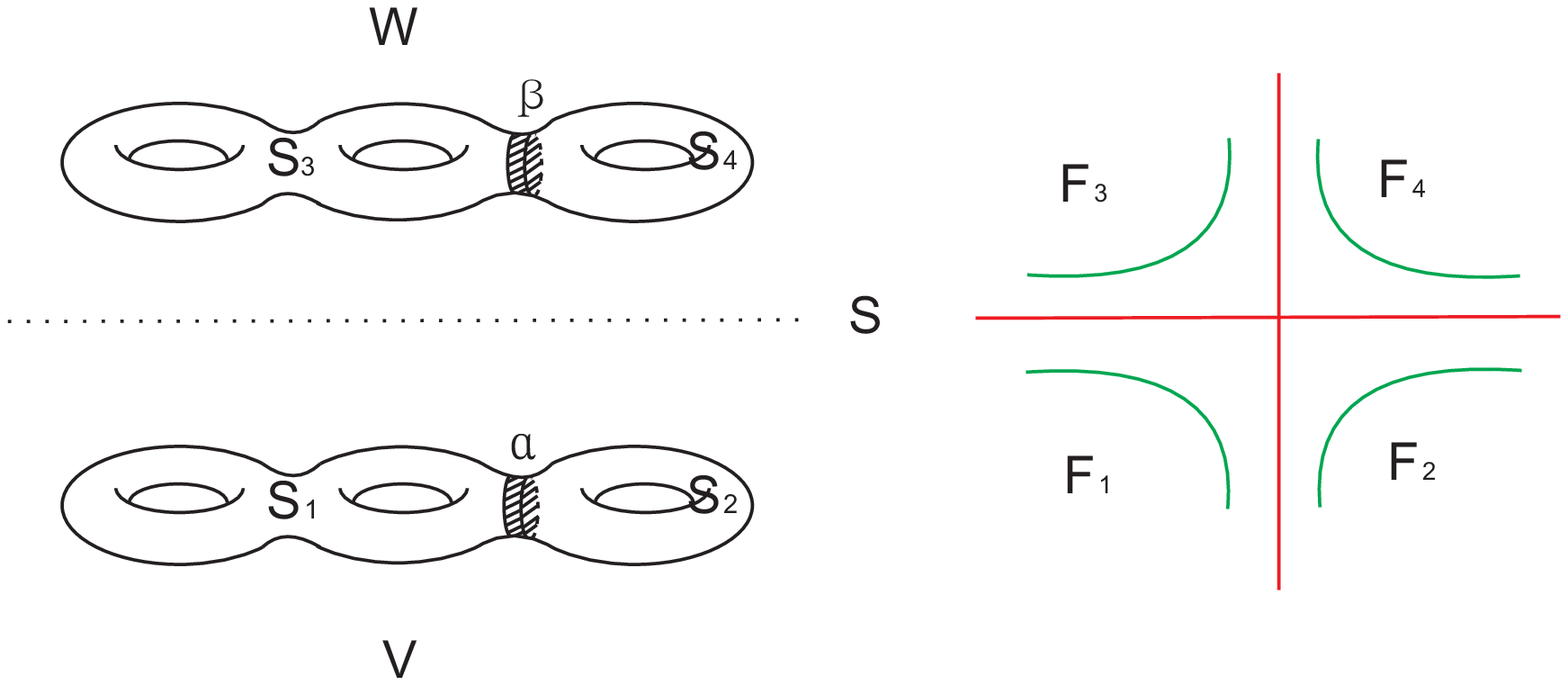}
\begin{center}
Figure 2
\end{center}
\end{center}

Let $F_{1}$ and $F_{2}$ be the  components of $\partial_{-} V$, and $S_{1}$ and $S_{2}$ be the two components of $S-\alpha$. Similarly, let $F_{3}$ and $F_{4}$ be the components of $\partial_{-} W$, and $S_{3}$ and $S_{4}$ be the two components of $S-\beta$. Now $B$ cuts $V$ into two manifolds $F_{1}\times I$ and $F_{2}\times I$, and $D$ cuts $W$ into two manifolds $F_{3}\times I$ and $F_{4}\times I$. See Figure 2. By Lemma 2.6, we may assume that $S_{3}$ is a once-punctured torus.

We first consider the compression body $V$. We may assume that $F_{i}=F_{i}\times\bigl\{0\bigr\}$, $S_{i}\cup B=F_{i}\times \bigl\{1\bigr\}$ for $1\leq i\leq 2$. Let $f_{F_{i}}: S_{i}\cup B\rightarrow F_{i}$ be the natural homeomorphism such that $f_{F_{i}}(x\times\bigl\{1\bigr\})=x\times\bigl\{0\bigr\}$ for $i=1, 2$ and $f_{F_{i}}(\emptyset)=\emptyset$. No doubt that $f_{F_{i}}$ is well defined. Then, for any two essential simple closed curves $\zeta, \theta\subset S_{i}\cup B$, $d_{\mathcal {C}(F_{i})}(f(\zeta),f(\theta))= d_{\mathcal {C}(S_{i}\cup B)}(\zeta, \theta)$ for $i=1, 2$. See Figure 3. Hence $f_{F_{i}}$ induces an isomorphism from $\mathcal {C}(S_{i}\cup B)$ to $\mathcal {C}(F_{i})$, for any $i=1,2$. Denote the isomorphism by $f_{F_{i}}$ too. Note that the shadow disk in Figure 3 is $B$.

\begin{center}
\includegraphics[height=2.5cm, width=7cm]{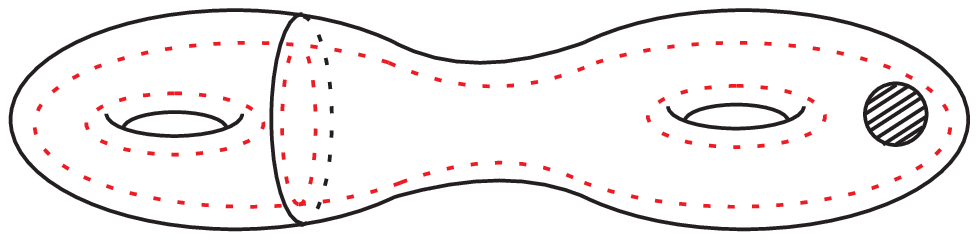}
\begin{center}
Figure 3
\end{center}
\end{center}\vskip 3mm

Let $\iota: S_{i}\rightarrow S_{i}\cup B$ be the inclusion map for $i=1, 2$. Note that $\partial S_{i}$ contains only one component. If $c$ is an essential simple closed curve in $S_{i}$, $\iota(c)$ is also essential in $S_{i}\cup B$.
Now, for any two essential simple closed curves $\zeta,\theta\subset S_{i}$
$d_{\mathcal {C}(S_{i}\cup B)}(\iota(\zeta), \iota(\theta))\leq d_{S_{i}}(\zeta,\theta)$ for $i=1, 2$.
Hence $\iota$ induces a distance non-increasing map from $\mathcal {C}(S_{i})$ to $\mathcal {C}(S_{i}\cup B)$, for any $i=1,2$. Denote the inclusion map by $\iota$ too. Then we can define a projection map :

\begin{center}
$ \psi_{F_{i}}=f_{F_{i}}\circ \iota\circ \pi_{S_{i}}: \mathcal {C}^{0}(S)\rightarrow \mathcal {C}^{0}(F_{i})$.
\end{center}

Since $d_{\mathcal {C(S)}}(\alpha, \beta)=n\geq 2$, $\alpha\cap \beta\neq \emptyset$. By the argument
in Section 2, $diam_{S_{i}}(\pi_{S_{i}}(\beta))\leq 2$. Hence $diam_{\mathcal {F}_{i}}(\psi_{F_{i}}(\beta))\leq 2$. \vskip 2mm

We start to attach a handlebody to $V$ along $F_{1}$. Then either

(1) $F_{1}$ is a torus. By Lemma 2.1, there is an essential simple closed curve $r$ in $F_{1}$ such that  $d_{\mathcal {C}(F_{1})}(\psi_{F_{1}}(\beta), r)\geq \mathcal {M}+1$. let $J_{r}$ be a solid torus such that $\partial J_{r}=F_{1}$, and $r$ bounds a disk in $J_{r}$. In this case, $J_{r}$ contains only one essential disk up to isotopy. Let $V_{F_{1}}$ be the manifold $V\cup J_{r}$. Or,

(2) $g(F_{1})\geq 2$. By Lemma 2.5, there is a full simplex $X$ on of $\mathcal {C}(F_{1})$ such that $d_{\mathcal {C}(F_{1})}(\mathcal {D}(H_{X}), \psi_{F_{1}}(\beta))\geq \mathcal {M}+1$, where $H_{X}$ is the handlebody obtained by attaching 2-handles to $F_{1}$ along $X$ then 3-handles to cap off the possible 2-spheres. In this case, we denote by $V_{F_{1}}$ the manifold $V\cup H_{X}$. \vskip 2mm

In whole words, $V_{F_{1}}$ is a compression body with only one minus boundary component $F_{2}$. See Figure 4. Hence $V_{F_{1}}\cup_{S} W$ is a Heegaard splitting. \vskip 2mm

\begin{center}
\includegraphics[height=4cm, width=6cm]{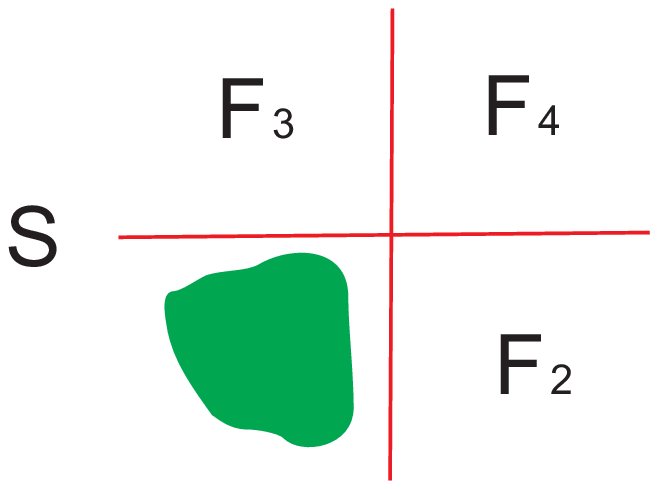}
\begin{center}
Figure 4
\end{center}
\end{center}\vskip 3mm

{\bf Claim 3.2.} The Heegaard distance of $V_{F_{1}}\cup_{S} W$, say $d_{\mathcal {C}(S)}(V_{F_{1}},W)$, is $n$. \vskip 2mm

{\bf Proof.} Suppose, otherwise, that $d_{\mathcal {C}(S)}(V_{F_{1}},W)=k<n$. Since $W$ contains only one essential disk $D$ up to isotopy such that $\partial D=\beta$,  there is an essential disk $B_{1}$ in $V_{F_{1}}$ such that
$d_{\mathcal {C}(S)}(\partial B_{1}, \beta)=k\leq n-1$, i.e, there is a geodesic $\mathcal {G}=\{a_{0}=\beta,...,a_{k}=\partial B_{1}\}$,
where $k\leq n-1$. \vskip 2mm

{\bf Fact 3.3.} $a_{j}\cap S_{1}\neq \emptyset$, for any $0\leq j\leq k$.\vskip2mm

Suppose that $a_{j}\cap S_{1}=\emptyset$ for some $0\leq j\leq k$. $j\neq k$ since $a_{k}=\partial B_{1}$, and $\partial S_{1}=\alpha$ and if $a_{k}\cap S_{1}=\emptyset$, then $B_{1}\subset F_{2}\times I$, and $B_{1}$ is inessential in $V_{F_{1}}$. $j\neq 0$ since $a_{0}=\beta$. Hence there is a geodesic $\mathcal {G}^{*}=\{\beta=a_{0},...,a_{j}, \alpha\}$. It means that $d_{\mathcal {C}(S)}(\alpha,\beta)\leq k<n$, a contradiction. \vskip 2mm

By Lemma 2.3,  $d_{\mathcal {C}(S_{1}\cup B)}(\partial B_{1}, \beta)\leq \mathcal{M}$. Furthermore, $d_{\mathcal {C}(F_{1})}(\psi_{F_{1}}(\partial B_{1}),\psi_{F_{1}}(\beta))\leq \mathcal{M}$. Depending on the way of intersection between $B_{1}$ and $B$, either \vskip 2mm

(1) $B_{1}\cap B=\emptyset$.  Since $B_{1}$ is not isotopic to $B$, $\psi_{F_{1}}(\partial B_{1})$ bounds an essential disk in $H_{X}$ or $J_{r}$ depending on $g(F_{1})$, where $H_{X}$ and $J_{r}$ are constructed as above. It contradicts the choice of $X$ or $r$. Or,\vskip 2mm

(2) $B_{1}\cap B \neq \emptyset$. Let $a$ be an outermost arc of $B_{1}\cap B$ on $B_{1}$. It means that $a$, together with a sub-arc $\gamma\subset\partial B_{1}$, bounds a disk $B_{\gamma}$ such that $B_{\gamma}\cap B=a$. Since $B$ cuts $V_{F_{1}}$ into a handlebody $H$ which  contains $F_{1}$ and a $I$-bundle $F_{2}\times I$, $B_{\gamma}\subset H$. Hence $\psi_{F_{1}}(\partial B_{1})$ bounds an essential disk in $H_{X}$ or $J_{r}$. By the argument in (1), it is impossible.  END.(Claim 1)

\begin{center}
\includegraphics[height=4cm, width=6cm]{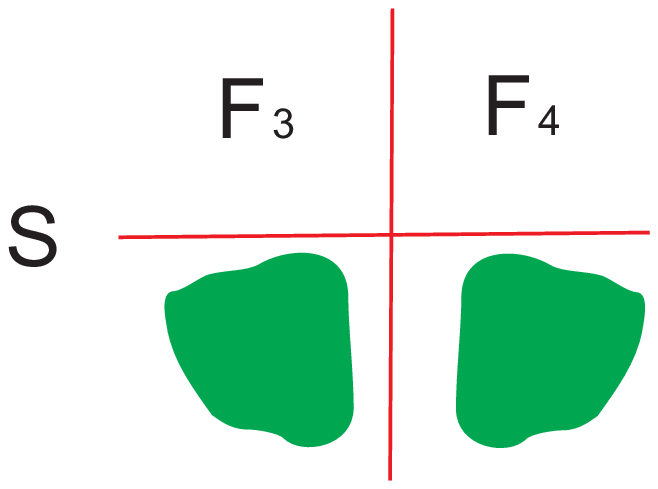}
\begin{center}
Figure 5
\end{center}
\end{center}\vskip 3mm

Now $V_{F_{1}}$ is a compression body which has only one minus boundary component $F_{2}$. Since $d_{\mathcal {C(S)}}(\alpha, \beta)=n\geq 2$, $\beta\cap S_{2}\neq\emptyset$. By Lemma 2.6, there is always a full simplex $Y$ on $F_{2}$ such that $d_{\mathcal {C}(F_{2})}(\mathcal {D}(H_{Y}), \psi_{F_{2}}(\beta))\geq \mathcal{M}+1$, where $H_{Y}$ is the handlebody obtained by attaching 2-handles to $F_{2}$ along $Y$ then 3-handles to cap off the possible 2-spheres, and $\psi_{F_{2}}$ is defined as before. Let $V_{F_{1},F_{2}}$ be the manifold obtained by  attaching  $H_{Y}$ to $V_{F_{1}}$ along $F_{2}$. See Figure 5. Then $V_{F_{1},F_{2}}$ is a handlebody. Hence $V_{F_{1},F_{2}}\cup_{S} W$ is also a Heegaard splitting. \vskip 2mm

{\bf Claim 3.4.} The Heegaard distance of $V_{F_{1},F_{2}}\cup_{S} W$, said $d_{\mathcal {C}(S)}(V_{F_{1},F_{2}},W)$, is $n$. \vskip 2mm

{\bf Proof.}  Suppose, otherwise, that $d_{\mathcal {C}(S)}(V_{F_{1},F_{2}},W)=k<n$. Since $W$ contains only one essential disk $D$ up to isotopy such that $\partial D=\beta$,
there is an essential disk $B_{2}$ in $V_{F_{1},F_{2}}$ such that $d_{\mathcal {C}(S)}(\partial B_{2}, \beta)=k$,
i.e, there is a geodesic $\mathcal {G}=\{a_{0}=\beta,...,a_{k}=\partial B_{2}\}$, where $k\leq n-1$. By the proof of Claim 1, $a_{j}\cap S_{2}\neq\emptyset$ for $0\leq j\leq k-1$. \vskip 3mm

 Note that $\partial B=\alpha$. Depending on the way of intersection between $B_{2}$ and $B$, either \vskip 2mm

(1) $B_{2}\cap B=\emptyset$. Since $d_{\mathcal {C}(S)}(\alpha,\beta)=n>k$, $B_{2}$ is not isotopic to $B$. By the proof of Claim 1, $\partial B_{2}$ does not lie in $S_{1}$.
Hence $\partial B_{2}\subset S_{2}$. It implies that $\psi_{F_{2}}(\partial B_{2})$ bounds an essential disk in $H_{Y}$. By lemma 2.3, $d_{\mathcal{C}(S_{2})}(\partial B_{2}, \beta)\leq \mathcal{M}$. Hence $d_{\mathcal{C}(F_{2})}(\psi_{F_{2}}(\partial B_{2}),\psi_{F_{2}}(\beta))\leq \mathcal{M}$,
and $d_{\mathcal {C}(F_{2})}(\mathcal{D}(H_{Y}), \psi_{F_{2}}(\beta))\leq \mathcal {M}$.
It contradicts the choice of $Y$. Or,\vskip3mm

(2) $B_{2}\cap B \neq \emptyset$. Let $a^{*}$ be an outermost arc of $B_{2}\cap B$ on $B_{2}$. This means that $a^{*}$, together with a sub-arc $\gamma^{*}\subset \partial B_{2}$, bounds a disk $B_{\gamma^{*}}$ such that $B_{\gamma^{*}}\cap B=a^{*}$.
By the proof of Claim 4.1 (2), $\gamma^{*}\subset S_{2}$. Thus $\psi_{F_{2}}(\partial B_{2})$ bounds an essential disk in $H_{Y}$. By the same argument in Claim 1, it is impossible. END. (Claim 2) \vskip3mm

Until now, we get a distance $n$ genus $g$ Heegaard splitting $V_{F_{1},F_{2}}\cup_{S} W$. In this case, $V_{F_{1},F_{2}}$ is a handlebody, and $W$ contains only one essential disk $D$ such that $\partial D=\beta$. Furthermore, we can cut $S$ along $\beta$ into two components $S_{3}$ and $S_{4}$, and cut $W$ along $D$ into two manifolds $F_{3}\times I$ and $F_{4}\times I$ such that $F_{i}=F_{i}\times\bigl\{0\bigr\}$, and $S_{i}\cup D=F_{i}\times\bigl\{1\bigr\}$ for $i=3, 4$. Now the shadow disk in Figure 3 is $D$. Let $f_{F_{i}}: S_{i}\cup D\rightarrow F_{i}$ be the natural homeomorphism such that  $f_{F_{i}}(x\times\bigl\{1\bigr\})=x\times\bigl\{0\bigr\}$ for $i=3, 4$. Then, for any two essential simple closed curves $\zeta, \theta\subset S_{i}\cup D$, $d_{\mathcal {C}(F_{i})}(f(\zeta),f(\theta))= d_{\mathcal {C}(S_{i}\cup D)}(\zeta, \theta)$ for $i=3, 4$, see Figure 3.
Hence $f_{F_{i}}$ induces an isomorphism from $\mathcal {C}(S_{i}\cup B)$ to $\mathcal {C}(F_{i})$, for any $i=1,2$. Denote the isomorphism by $f_{F_{i}}$ too.

Let $\iota: S_{i}\rightarrow S_{i}\cup D$ be the inclusion map for $i=3, 4$. Note that $\partial S_{i}$ contains only one component. If $c$ is an essential simple closed curve in $S_{i}$, $\iota(c)$ is also essential in $S_{i}\cup D$. Now, for any two essential simple closed curves $\zeta,\theta\subset S_{i}$
$d_{\mathcal {C}(S_{i}\cup D)}(\iota(\zeta), \iota(\theta))\leq d_{S_{i}}(\zeta,\theta)$ for $i=3, 4$.
Hence $\iota$ induces a distance non-increasing map from $\mathcal {C}(S_{i})$ to $\mathcal {C}(S_{i}\cup B)$, for any $i=1,2$. Denote the inclusion map by $\iota$ too. Then we can define a projection map :

\begin{center}

$ \psi_{F_{i}}=f_{F_{i}}\circ \iota\circ \pi_{S_{i}}: \mathcal {C}^{0}(S)\rightarrow \mathcal {C}^{0}(F_{i})$. \end{center}
\vskip 3mm

Since $n\neq 2$, there are two cases: \vskip 2mm

{\bf Case 1.}  $n\geq 3$. \vskip 2mm

Since $V_{F_{1},F_{2}}\cup_{S} W$  is a distance $n$ Heegaard splitting of  genus $g$, and $W$ contains only an essential disk $D$ up to isotopy, $S_{3}$ and $S_{4}$ are incompressible in $V_{F_{1},F_{2}}$. Hence $\beta=\partial S_{3}=\partial S_{4}$ is disk-busting in $V_{F_{1},F_{2}}$. Since $g\geq 3$, and $g(S_{3})=1$, $V_{F_{1},F_{2}}$ is
not an I-bundle over some compact surface with  $S_{i}$  a horizontal boundary of the I-bundle,
and the vertical boundary of this I-bundle  a single annulus for $i=3, 4$. By Lemma 2.4, $diam_{S_{i}}(\mathcal {D}(V_{F_{1},F_{2}}))\leq 12$ for $i=3, 4$. Hence $diam_{F_{i}}(\psi_{F_{i}}(\mathcal {D}(V_{F_{1},F_{2}})))\leq 12$. \vskip 3mm

Since $F_{3}$ is a torus, by Lemma 2.1, there is an essential simple closed curve $\delta$ in $F_{3}$ such that  $d_{\mathcal {C}(F_{3})}(\psi_{F_{3}}(\mathcal {D}(V_{F_{1},F_{2}})), \delta)\geq \mathcal {M}+1$. Let $W_{F_{3}}$ the be the manifold obtained attaching a solid $J_{\delta}$ to $W$ along $F_{3}$ so that $\delta$ bounds a disk in $J_{\delta}$. Then $W_{F_{3}}$ is a compression body.

Since $g\geq 3$, $g(F_{4})\geq 2$. By Lemma 2.5, there is a full simplex $Z$ of $\mathcal {C}(F_{4})$ such that $d_{\mathcal {C}(F_{4})}(\mathcal {D}(H_{Z}), \psi_{F_{4}}(\mathcal {D}(V_{F_{1},F_{2}})))\geq \mathcal {M}+1$, where $H_{Z}$ is the handlebody obtained by attaching 2-handles to $F_{4}$ along $Z$ then 3-handles to cap off the possible 2-spheres. In this case, let $W_{F_{3}, F_{4}}$ be the handlebody $W_{F_{3}}\cup H_{Z}$. Now $V_{F_{1}, F_{2}}\cup_{S} W_{F_{3}, F_{4}}$ is a Heegaard splitting of a closed 3-manifold.\vskip 2mm

\begin{center}
\includegraphics[height=6cm, width=8cm]{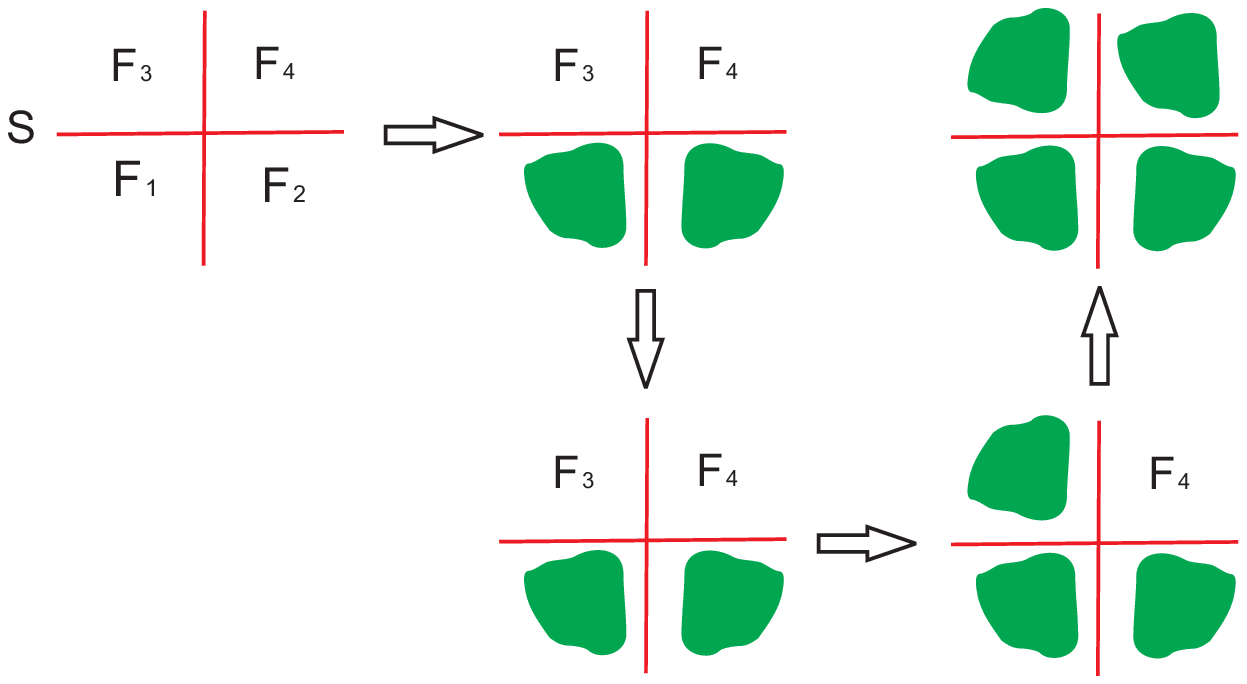}
\begin{center}
Figure 6
\end{center}
\end{center}

{\bf Claim 3.5.} The distance of $V_{F_{1}, F_{2}}\cup_{S} W_{F_{3}, F_{4}}$, said $d_{\mathcal {C}(S)}(V_{F_{1},F_{2}},W_{F_{3},F_{4}})$, is $n$. \vskip 2mm

{\bf Proof.}Let $D$ be the essential disk in $W_{F_{3},F_{4}}$ bounded by $\beta$. Suppose, otherwise, that $d=k<n$. Then there is a geodesic  $\mathcal{G}=\{a_{0}=\partial B_{1},...,a_{k}=\partial D_{1}\}$, where $k\leq n-1$, $B_{1}$ is a disk in $V_{F_{1}, F_{2}}$, and $D_{1}$ is a disk in $W_{F_{3}, F_{3}}$. $\alpha_{i}\cap \beta\neq \emptyset$, for any $0\leq i\leq k-1$ for if not, the distance of $V_{F_{1}, F_{2}}\cup_{S} W$ would be at most $k<n$. Similarly, $D_{1}$ is not isotopic to $D$. \vskip 3mm

Then either

(1) $D_{1}\cap D=\emptyset$. Then $\partial D_{1}$ lies in one of $S_{3}$ and $S_{4}$, say $S_{3}$. Hence $\psi_{F_{3}}(\partial D_{1})=\delta$. By Lemma 2.3, $diam_{S_{3}}(\mathcal {D}(\mathcal{G}))\leq \mathcal {M}$. Since $\pi_{S_{3}}(\partial B_{1})\in \pi_{S_{3}}(\mathcal{D}(V_{F_{1}, F_{2}}))$, \\$d_{\mathcal{C}(S_{3})} (\pi_{S_{3}}(\mathcal{D}(V_{F_{1}, F_{2}})), \partial D_{1})\leq \mathcal{M}$. Hence $d_{\mathcal {C}(F_{3})}(\psi_{F_{3}}(\mathcal {D}(V_{F_{1},F_{2}})), \psi_{F_{3}}(\partial D_{1})=\delta)\leq \mathcal {M}$, a contradiction.  Or,\vskip3mm

(2) $D_{1}\cap D \neq \emptyset$, Let $c$ be an outermost arc of $B_{2}\cap B$ on $B_{2}$. This means that $c$, together with a sub-arc $\delta^{*}\subset \partial B_{2}$, bounds a disk $D_{c}$ such that $D_{c}\cap D=\gamma^{*}$. We may assume that $\partial D_{c}\subset S_{4}$. By Lemma 2.3, $diam_{S_{4}}(\mathcal {D}(\mathcal{G}))\leq \mathcal {M}$. Hence $d_{\mathcal {C}(F_{4})}(\psi_{F_{4}}(\mathcal {D}(V_{F_{1},F_{2}})), \psi_{F_{4}}(\partial B_{2}))\leq \mathcal {M}$. Note that $\psi_{F_{4}}(\partial B_{2})\in \mathcal {D}(H_{Z})$. By the same argument in (1), it is impossible. END.  \vskip3mm

Now we suppose that $n=1$. \vskip 2mm

Let $M_{1}$ and $M_{2}$ be two 3-manifolds with homeomorphic connected boundary. Let $M^{f}$ be the manifold obtained by gluing $M_{1}$ and $M_{2}$ along a homeomorphism from $\partial M_{1}$ to $\partial M_{2}$. Let $M_{i}=V_{i}\cup_{S_{i}} W_{i}$ be a minimal Heegaard splitting for $i=1, 2$. In this case, $M^{f}$ has a natural Heegaard called the amalgamation of $V_{1}\cup_{S_{1}} W_{1}$ and $V_{2}\cup_{S_{2}} W_{2}$. The following facts are well known: \vskip 2mm

(1) If the gluing map $f$ is enough complicated, then  the amalgamation of $V_{1}\cup_{S_{1}} W_{1}$ and $V_{2}\cup_{S_{2}} W_{2}$ is unstabilized, see \cite{bss}, \cite{La02}, \cite{l05}, \cite{s}.

(2) If both $V_{1}\cup_{S_{1}} W_{1}$ and $V_{2}\cup_{S_{2}} W_{2}$ have high distance, then the amalgamation of $V_{1}\cup_{S_{1}} W_{1}$ and $V_{2}\cup_{S_{2}} W_{2}$ is unstabilized, See \cite{KQ}, \cite{YL}.

Now let $M_{i}=V_{i}\cup_{S_{i}} W_{i}$ be a Heegaard splitting of genus two such that $\partial M_{i}$ is a torus, and $d(S_{i})>8$ for $i=1, 2$, then, by the main result in \cite{KQ}, the amalgamation  of $V_{1}\cup_{S_{1}} W_{1}$ and $V_{2}\cup_{S_{2}} W_{2}$, say $V\cup_{S} W$, is unstabilized. Furthermore, $g(S)=3$. \vskip 2mm

Suppose that $g\geq 4$. By the above argument, there are a Heegaard splitting $M_{1}=V_{1}\cup_{S_{1}} W_{1}$ of genus $g-1$ such that $g(\partial M_{1})=2$, and $d(S_{1})\geq 2g$, and a Heegaard splitting $V_{2}\cup_{S_{2}} W_{2}$ of genus 3 such that $g(\partial M_{2})=2$, and $d(S_{2})\geq 2g$. Hence both $M_{1}$ and $M_{2}$ are hyperbolic. By  the main result in \cite{KQ}, the amalgamation  of $V_{1}\cup_{S_{1}} W_{1}$ and $V_{2}\cup_{S_{2}} W_{2}$, say $M=V\cup_{S} W$, is unstabilized.
Furthermore, $g(S)=g$. By Thurston's Theorem, $M$ is hyperbolic.
\vskip 3mm

END(Proposition 3.1)\vskip 5mm

{\bf Remark.} The strongly irreducible Heegaard splitting $V\cup_{S} W$  where both $V$ and $W$ contain only one essential separating disk up to isotopy independently is always a minimal Heegaard splitting of $M=V\cup_{S}W$. T.Li \cite{l05} defined a sub-complex $\mathcal {U}(F_{1})$, for $F_{1}\subset \partial _{-}V$ and proved that for any handlebody $H$ attached to $M$ along $F_{1}$, if $d_{\mathcal {C}(F_{1})}(\mathcal {U}(F_{1}), \mathcal {D}(H))$ is larger than a constant $\mathcal {K}$ which depends on $M $ and $H$, then the new generated Heegaard splitting $V_{F_{1}}\cup_{S} W$ is still the minimal Heegaard splitting of $M^{F_{1}}=V_{F_{1}}\cup_{S} W$. Similar to the other boundaries of $M$.
Now in our construction of distance $n\geq 2$ strongly irreducible Heegaard splitting (for n=2, see section 5), we can choose a full simplex $X$ in $F_{1}$ such that $d_{\mathcal {C}(F_{1})}(\psi_{F_{1}}(\mathcal {D}(W)), \mathcal {D}(H_{X}))$ is large enough and $d_{\mathcal {C}(F_{1})}(\mathcal {U}(F_{1}), \mathcal {D}(H_{X}))$ is larger than $\mathcal {K}$. Then the new Heegaard splitting $V_{F_{1}}\cup_{S} W$ is still the minimal Heegaard splitting of $M^{F_{1}}=V_{F_{1}}\cup_{S} W$ and has the same distance as the older one.

\section{Proof of Theorem 2}
We will prove Theorem 2 in this section.\vskip 3mm

{\bf Theorem 2.} For any integers $g\geq 2$ and $n\geq 4$, there are infinitely many non-homeomorphic closed 3-manifolds which admit distance $n$ Heegaard splittings of genus $g$. \vskip 2mm

{\bf Proof.} Let $S_{g}$ be a closed surface of genus $g$. By Lemma 2.6, for each $m\geq 2$,  there is a geodesic $\mathcal {G}^{m}=\{\alpha=a_{0}^{m}, a_{1}^{m}, ..., a_{n-1}^{m}, a_{n}^{m}=\beta^{m}\}$ in $\mathcal {C}(S_{g})$ such that

(1)  $a_{i}$ is non-separating in $S_{g}$ for  $1\leq i\leq n-1$,  $\alpha$ and $\beta^{m}$ are two essential separating simple closed curves on $S_{g}$ for $m\geq 2$,

(2) $m\mathcal{M}+2\leq d_{\mathcal {C}(S^{a_{i}})}(a_{i-1}, a_{i+1})= m\mathcal{M}+6$,
where $S^{a_{i}}$ is the surface $S-N(a_{i})$ for $ 1\leq i \leq n-1$, and

(3) one component of $S_{g}-\beta^{m}$ has genus one. \vskip 2mm

Without loss of generality, we assume that $\mathcal{M}\geq 6$. Let $M_{m}$ be the manifold obtained by attaching two 2-handles to $S_{g}\times [-1,1]$ along $\alpha\times\{-1\}$ and $\beta^{m}\times\{1\}$. We denote also by $S_{g}$ the surface $S_{g}\times\{0\}$. Now $M_{m}$ has a Heegaard splitting as $V_{m}\cup_{S_{g}} W_{m}$, where $V_{m}$ is the compression body obtained by attaching a 2-handle to $S\times [-1,0]$ along $\alpha\times\{-1\}$, and $W_{m}$ is the manifold obtained by attaching a 2-handle to $S\times [0,1]$ along $\beta^{m}\times\{1\}$. Then $\partial_{-}V_{m}$ contains two components $F_{1}$ and $F_{2}$, and $\partial_{-}W_{m}$ contains two components $F_{3}^{m}$ and $F_{4}^{m}$. See Figure 7. Furthermore, one component of $F_{3}^{m}$ and $F_{4}^{m}$ has genus one.

\begin{center}
\includegraphics[height=3cm, width=6cm]{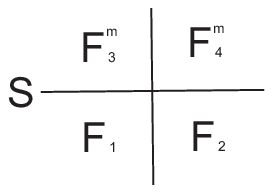}
\begin{center}
Figure 7
\end{center}
\end{center}

By the proof of Theorem 1(1), there is a closed 3-manifold $M_{m}^{*}$ which admits a distance $n$ Heegaard splitting $V_{m}^{*}\cup_{S_{g}} W_{m}^{*}$, where $V_{m}^{*}$ is obtained by attaching handlebodies $H_{X_{1}}$ and $H_{X_{2}}$ to $V_{m}$ along $F_{1}$ and $F_{2}$, and $W_{m}^{*}$ is obtained by attaching handlebodies $H_{Y_{1}}$ and $H_{Y_{2}}$ to $W_{m}$ along $F_{3}^{m}$ and $F_{4}^{m}$ such that

(1) $d_{\mathcal{C}(F_{i})}(\psi_{F_{i}}(\beta^{m}), \mathcal{D}(H_{X_{i}}))\geq \mathcal{M}+15$ for $i=1, 2$, and

(2) $d_{\mathcal{C}(F_{i})}(\psi_{F_{i}}(\alpha), \mathcal{D}(H_{Y_{i}}))\geq \mathcal{M}+15$ for $i=3, 4$. \vskip 3mm

Replace $M_{m}^{*}$, $V_{m}^{*}$ and $W_{m}^{*}$ by $M_{m}$, $V_{m}$ and $W_{m}$.
Now $\mathcal {G}^{m}=\{\alpha=a_{0}^{m}, a_{1}^{m},..., a_{n-1}^{m}, \\ a_{n}^{m}=\beta^{m}\}$ is also a geodesic of $\mathcal{C}(S_{g})$ realizing the distance of $M_{m}=V_{m}\cup_{S_{g}} W_{m}$. \vskip 3mm

%By the proof of theorem 1 part I, for
%any $\mathcal {G}^{t}\in B$, we can attaching suitable Handlebodies and get a distance $n$ Heegaard splitting $V^{t}\cup_{S} W^{t}$. Hence

%We call a
%geodesic $\mathcal {G}=\{a_{0},...,a_{n}\}$ realizing the Heegaard distance of a Heegaard splitting $V\cup_{S} W$ if
%the Heegaard distance is $n$, $a_{0}$(resp. $a_{n}$) bounds an essential disk in $V$ (resp. $W$).\vskip3mm

{\bf Claim 4.2.} Let $\mathcal {G}=\{b_{0},...,b_{n}\}$ be a geodesic of $\mathcal{C}(S_{g})$ realizing the distance of $V_{m}\cup_{S_{g}} W_{m}$. Then $b_{i}= a_{i}^{m}$ for any $1\leq i \leq n-1$. \vskip 2mm

{\bf Proof.} Let $S_{1}$ and $S_{2}$ be the two components of $S_{g}-\alpha$. We assume that $b_{0}$ bounds a disk $B_{0}$ in $V_{m}$, and $b_{n}$ bounds a disk $D_{n}$ in $W_{m}$. We first prove that $\alpha$(resp. $\beta^{m}$) is disjoint from $b_{1}$ (resp. $b_{n-1}$).\vskip3mm

Suppose, otherwise, that $\alpha\cap b_{1}\neq\emptyset$. Hence $b_{0}$ is not isotopic to $a_{0}^{m}=\alpha$.
Then either

(1) $B_{0}\cap B\neq \emptyset$. Let $a$ be an outermost arc of $B_{0}\cap B$ on $B_{0}$. It means that $a$, together a sub-arc of $\gamma\subset \partial B_{0}$, bounds a disk $B_{\gamma}$ such that $B_{\gamma}\cap B=a$. Without assumption, we may assume that $\gamma\subset S_{1}$. By the argument in section 3, $\psi_{F_{1}}(\partial B_{0})$ bounds an essential disk in $H_{X_{1}}$. But with $b_{1}\cap \partial S_{1}\neq \emptyset$, it implies that
$d_{\mathcal {C}(S_{1})}(b_{0},b_{n})\leq \mathcal {M}$. Hence $d_{\mathcal {C}(F_{1})}(\psi_{F_{1}}(b_{n}), \mathcal {D}(H_{X_{1}}))\leq \mathcal {M}$. Or,

(2) $B_{0}\cap B=\emptyset$. By $b_{1}\cap \alpha\neq \emptyset$, $B_{0}$ is not isotopic to $B$. Then
$\partial B_{0}$ is essential in $S_{1}$ or $S_{2}$. We assume that $\partial B_{0}\subset S_{1}$. The other case is
similar. Hence by (1), $d_{\mathcal {C}(F_{1})}(\psi_{F_{1}}(b_{n}), \mathcal {D}(H_{X_{1}}))\leq \mathcal {M}$.

However, by Heegaard distance is at least 4 and $\alpha=\partial S_{1}=\partial S_{2}$ bounds an essential disk in $V^{m}$, it means that $\alpha$ is disk-busting for $W^{m}$ and $W^{m}$ can not be the I-bundle of compact surface with $S_{1}$ or $S_{2}$ as one of its horizontal boundary. Then by Lemma 2.4, $diam_{\mathcal {C}(S_{1})}(\mathcal {D}(W^{m}))\leq 12$ and $diam_{\mathcal {C}(S_{2})}(\mathcal {D}(W^{m}))\leq 12$. Hence $diam_{\mathcal {C}(F_{1})}(\mathcal {D}(W^{m}))\leq 12$ and $diam_{\mathcal {C}(F_{2})}(\mathcal {D}(W^{m}))\leq 12$. Together with (1) and (2), by triangle inequality, $d_{\mathcal{C}(F_{1})}(\psi_{F_{1}}(\beta^{m}), \mathcal{D}(H_{X_{1}}))\leq \mathcal{M}+12$.
It contradicts the choice of $X_{1}$ in $F_{1}$. The other case is similar.
\vskip 3mm

Let $\mathcal {G}^{*}=\{\alpha=a_{0}^{m}, b_{1},...,b_{n-1}, a_{n}^{m}\}$ be a new geodesic realizing the distance of $V_{m}\cup_{S_{g}}W_{m}$. Now we prove that $b_{1}$ is isotopic to
$a_{1}^{m}$.

Suppose, otherwise, that $b_{1}$ is not isotopic to
$a_{1}^{m}$. Note that $b_{i}$ is not isotopic to $a_{1}^{m}$. Otherwise, the distance of  $V_{m}^{c}\cup_{S_{g}}W_{m}^{c}$ would be at most $n-1$. Let $S^{a_{1}^{m}}$ be the surface $S_{g}-N(a_{1}^{m})$, where $N(a_{1}^{m})$ is a open regular neighborhood of $a_{1}^{m}$ on $S_{g}$. By Lemma 2.3, $d_{\mathcal{C}(S^{a_{1}^{m}})}(\pi_{S^{a_{1}^{m}}}(a_{0}^{m}),\pi_{S^{a_{1}^{m}}}(a_{n}^{m}))\leq\mathcal{M}$. Now let's consider the shorter geodesic $\mathcal {G}^{**}=\{ a_{2}^{m}, ..., a_{n-1}^{m}, a_{n}^{m}=\beta^{m}\}$ which is a sub-geodesic of $\mathcal {G}^{m}=\{\alpha=a_{0}^{m}, a_{1}^{m}, ..., a_{n-1}^{m}, a_{n}^{m}=\beta^{m}\}$. Due to
the definition of geodesic in curve complex, $a_{i}^{m}$ is not isotopic to $a_{1}^{m}$ for any $i\geq 2$. By Lemma 2.3 again,  $d_{\mathcal{C}(S^{a_{1}^{m}})}(\pi_{S^{a_{1}^{m}}}(a_{2}^{m}),\pi_{S^{a_{1}^{m}}}(a_{n}^{m}))\leq\mathcal{M}$. Hence $d_{\mathcal{C}(S^{a_{1}^{m}})}(\pi_{S^{a_{1}^{m}}}(a_{0}^{m}),\pi_{S^{a_{1}^{m}}}(a_{2}^{m}))\leq 2\mathcal{M}$. This contradicts our assumption on $d_{\mathcal{C}(S^{a_{1}^{m}})}(\pi_{S^{a_{1}^{m}}}(a_{0}^{m}),\pi_{S^{a_{1}^{m}}}(a_{n}^{m}))$. Hence $b_{1}$ is isotopic to $a_{1}^{m}$.

By induction on $i$, the claim holds. End (Claim 4.2) \vskip 3mm

Replace $M_{m}=V_{m}\cup_{S_{g}} W_{m}$ by $M_{m}=V_{m}\cup_{S_{g}^{m}} W_{m}$.

The following claim reveals the connection between geodesics in curve complex and closed 3-manifolds. \vskip 2mm

{\bf Claim 4.3.} For any $2 \leq t\neq s\in N$, either

(1) $M_{t}=V_{t} \cup_{S_{g}^{t}} W_{t}$ and $M_{s}=V_{s}\cup_{S_{g}^{s}} W_{s}$ are two different 3-manifolds up to homeomorphism. Or,

(2) $M_{t}$ is homeomorphic to $M_{s}$, but $V_{t}\cup_{S_{g}^{t}} W_{t}$ and $V_{s}\cup_{S_{g}^{s}} W_{s}$ are two different Heegaard splittings of $M_{t}$ up to homeomorphic equivalence. \vskip 3mm

{\bf Proof.} Suppose that $M_{t}$ is homeomorphic to $M_{s}$ for some $2 \leq t\neq s\in N$. If (2) fails, then $V_{t}\cup_{S_{g}^{t}} W_{t}$ and $V_{s}\cup_{S_{g}^{s}} W_{s}$ are homeomorphic. It means that there is a homeomorphism $f$ from $M_{t}$ to $M_{s}$ such that $f((S_{g}^{t}; V_{t},W_{t}))=(S_{g}^{s}; V_{t},W_{t})$. We assume that
$f(V_{t})=V_{s}$ and $f(W_{t})=W_{s}$. The other case is similar. It is well known that $f$ induces an isomorphism from $\mathcal {C}(S_{g}^{t})$ to $\mathcal {C}(S_{g}^{s})$, still denoted by $f$. Then for the geodesic $\mathcal {G}^{t}=\{\alpha=a_{0}^{t}, a_{1}^{t}, ..., a_{n-1}^{t}, a_{n}^{t}=\beta^{t}\}$ which realizes the distance of $V_{t} \cup_{S_{g}^{t}} W_{t}$, $f(\mathcal{G})$ is also a geodesic in $\mathcal {C}(S_{g}^{s})$ realizing the distance of $V_{s}\cup_{S_{g}^{s}} W_{s}$. By Claim 4.2, $f(a_{j}^{t})$ is isotopic to $a_{j}^{s}$ for $1\leq j\leq n-1$.\vskip3mm

As $n\geq 4$, we choose $a^{t}_{2}$.
Since $f(a^{t}_{2})$ is isotopic to $a^{s}_{2}$, we can perform an isotopy on $S_{g}^{s}$ such that the composition of $f$ with the isotopy gives an homeomorphism $f^{\star}$ from $S_{t}$ to $S_{s}$ and $f^{\star}(a^{t}_{2})=a^{s}_{2}$. Even more, $f^{\star}(V_{t})=V_{s}$ and $f^{\star}(W_{t})=W_{s}$. It's still true that $f^{\star}$ induces an automorphism from $\mathcal {C}(S_{g}^{t})$ to $\mathcal {C}(S_{g}^{s})$, denoted by $f^{\star}$ too. Thus $f^{\star}(\mathcal {G}^{t})$
is also a geodesic realizing the distance of $V_{s}\cup_{S_{g}^{s}} W_{s}$. By Claim 4.2 again, for any $ 1\leq j\leq n-1$, $f^{\star}(a_{j}^{t})$ is still isotopic to $a_{j}^{s}$. Hence $f^{\star}(a_{1}^{t})$ (resp. $f^{\star}(a_{3}^{t})$) is isotopic to $a_{1}^{s}$ (resp. $a_{3}^{t}$).\vskip 3mm

Let $S^{a^{t}_{2}}$ be the surface $S_{g}^{t}-N(a^{t}_{2}))$, where $N(a^{t}_{2})$ is an open regular neighborhood of $a^{t}_{2}$ on $S_{g}^{t}$, and $S^{a^{s}_{2}}$ be the surface of $S_{g}^{s}-N(a^{s}_{2})$. Then $f^{\star} (S^{a^{t}_{2}})=S^{a^{s}_{2}}$ and $f^{\star}\mid_{S^{a^{t}_{2}}}$ is a homeomorphism. Hence $f^{\star}$ also induces an isomorphism from $\mathcal {C}(S^{a^{t}_{2}})$ to $\mathcal {C}(S^{a^{s}_{2}})$,
still denoted by $f^{\star}$. Now we can also assume $a_{1}^{t}\cap a_{2}^{t}=\emptyset$ and $a_{3}^{t}\cap a_{2}^{t}=\emptyset$. Thus $f^{\star}(a_{1}^{t})\cap (f^{\star}(a_{2}^{t})=a_{2}^{s})=\emptyset$ and $f^{\star}(a_{3}^{t})\cap (f^{\star}(a_{2}^{t})=a_{2}^{s}) =\emptyset$.
Then $d_{\mathcal{C}(S^{a^{t}_{2}})}(a_{1}^{t}, a_{3}^{t})=d_{\mathcal{C}(S^{a^{s}_{2}})}(f^{\star}(a_{1}^{t}), f^{\star}(a_{3}^{t}))$. On the other side, $f^{\star}(a_{1}^{t})$ (resp.$f^{\star}(a_{3}^{s})$) must be isotopic to $a_{1}^{s}$ (resp.$a_{3}^{s}$) in $S^{a^{s}_{2}}$ for
$\vartriangleright$ if not, then after removing possible Bigon capped by them, they bounds no annulus in $S^{a^{s}_{2}}$, thus they bounds no annulus and Bigon in $S_{g}^{s}$. By Bigon Criterion (proposition 1.7\cite{FM}), they realizes the geometry intersection number. Since they are isotopic in $S_{g}^{s}$, they must be disjoint in $S_{g}^{s}$. Hence they must bounds an annulus in $S_{g}^{s}$ $\vartriangleleft$.
So $d_{\mathcal{C}(S^{a^{t}_{2}})}(a_{1}^{t}, a_{3}^{t})=d_{\mathcal{C}(S^{a^{s}_{2}})}(f^{\star}(a_{1}^{t}), f^{\star}(a_{3}^{t}))=d_{\mathcal{C}(S^{a^{s}_{2}})}(a_{1}^{s}, a_{3}^{s})$. However, by the assumption [
$t\mathcal {M}+2\leq d_{\mathcal{C}(S^{a^{t}_{2}})}(a_{1}^{t}, a_{3}^{t})\leq t\mathcal {M}+6$, $s\mathcal {M}+2\leq d_{\mathcal{C}(S^{a^{s}_{2}})}(a_{1}^{t}, a_{3}^{t})\leq s\mathcal {M}+6$ and $\mathcal{M}\geq 6$ ],  $d_{\mathcal{C}(S^{a^{t}_{2}})}(a_{1}^{t}, a_{3}^{t})\neq d_{\mathcal{C}(S^{a^{s}_{2}})}(a_{1}^{s}, a_{3}^{s})$, a contradiction.
 End (Claim 4.3) \vskip 3mm

The Waldhausen conjecture proved by Johanson (\cite{j01},\cite{j02}) and  Li \cite{l02,l03} implies that, for any positive integer $g$, an atoroidal closed 3-manifold $M$ admits only finitely many Heegaard splittings of genus $g$ up to homeomorphism. Since $M_{t}$ admits a Heegaard splitting with distance at least 4, it is atoroidal for any $t\geq 2$, see \cite{ha} and \cite{sch01}. Now Theorem 2 is immediately from Claim 2 and the Waldhausen conjecture. END \vskip 3mm

\section{Proof of Theorem 1(2)}

We rewrite the second part  of Theorem 1 as the following proposition: \vskip 3mm

{\bf Proposition 5.1.} For any integer $g\geq 2$, there is a hyperbolic closed 3-manifold which admits a distance 2 Heegaard splitting of genus $g$. \vskip 3mm

{\bf Proof.} By the remark on Theorem 1, there is a hyperbolic closed 3-manifold which admits a distance 2 Heegaard splitting of genus 2.\vskip 2mm

%We first prove it for $g=2$. Let $K$ be a hyperbolic 2-bridge knot. Then its complement, say $E(K)$, admits a distance 2 Heegaard splitting of genus 2, say $V\cup_{S} W$. We may assume that $T=\partial E(K)\subset \partial_{-}W$.
%Then there are at most ten slopes $r$ on $T$ such that the manifold $E_{r}(K)$ obtained by doing Dehn filling on $E(K)$ along $r$ is non-hyperbolic, see \cite{th} \cite {A00}, \cite{A10}, \cite{La12}. So we can always choose a slope $c\subset T$ such that $E_{c}(K)$ is still hyperbolic. And $E_{c}(K)$ admits a genus 2 Heegaard splitting $V\cup_{S} W_{r}$. It is not hard to see that
%distance of $V\cup_{S} W_{r}$ is less than or equal to 2. We say distance of this Heegaard splitting can not be 0. Otherwise, it would be stabilized, and $E_{c}(K)$ should be a lens space. Distance of this Heegaard splitting
%can not be 1 since if it is, then this Heegaard splitting is weakly reducible. By A.Casson and C.Gordon's work \cite{CG}, and M.Scharlemann and A.Thompson's thin position \cite{st1}, this Heegaard splitting can be written as an amalgamation of some
%strongly irreducible Heegaard splittings. But a genus 2 Heegaard splitting can not be written as an amalgamation.\vskip 3mm

Suppose now that $g\geq 3$. \vskip 2mm

{\bf Assumption 1.} \ Let $S$ be a closed surface of genus $g$. By Lemma 2.6, there are two separating slopes $\alpha$ and $\gamma$ such that

(1)$d_{\mathcal {C}(S)}(\alpha, \gamma)=2$,

(2) one component of $S-\alpha$, say $S_{1}$,  has genus one while another component of $S-\alpha$, say $S_{2}$, has genus $g-1$,

(3)  one component of $S-\gamma$, say $S_{3}$, has genus one, while another component of $S-\gamma$, say $S_{4}$, has genus $g-1$

(4) there is a non-separating slope $\beta$ on $S$ such that
$\alpha$ and $\gamma$ are disjoint from $\beta$,  and $d_{\mathcal {C}(S^{\beta})}(\alpha, \gamma)>4$, where $S^{\beta}$ is the surface $S-\eta(\beta)$,  and

(5) $\beta \subset S_{2}\cap S_{4}$. \vskip 2mm

Let $V$ be the compression body obtained by attaching a separating 2-handle to $S\times [0,1]$  along $\alpha\times \bigl\{1\bigr\}$, and $W$ be the compression body obtained by attaching a separating 2-handle to $S\times [-1,0]$ along $\gamma\times \bigl\{-1\bigr\}$. Denote $S\times\bigl\{0\bigr\}$ by $S$ too. Then $V\cup_{S} W$ is a Heegaard splitting. Since $V$ contains only one essential disk $B$ with $\partial B=\alpha$ up to isotopy, and $W$ contains only one essential disk $D$ with $\partial D=\gamma$ up to isotopy, $d_{\mathcal {C}(S)}(V,W)=2$.

Let $F_{1}$ and $F_{2}$ be the  components of $\partial_{-} V$, such that $F_{i}$ is homeomorphic to $S_{i}\cup B$ for $i=1, 2$. Similarly, let $F_{3}$ and $F_{4}$ be the components of $\partial_{-} W$ such that $F_{i}$ is homeomorphic to $S_{i}\cup D$ for $i=3, 4$. Then both $S_{1}$ and $S_{3}$ are once-punctured tori, and $F_{1}$ and $F_{3}$ are two tori, see Figure 2. Furthermore, both $F_{3}$ and $F_{4}$ have genus at least 2. Now $B$ cuts $V$ into two manifolds $F_{1}\times I$ and $F_{2}\times I$, and $D$ cuts $W$ into two manifolds $F_{3}\times I$ and $F_{4}\times I$.

Since $d_{\mathcal {C}(S)}(V,W)=2$, $\gamma\cap S_{i}\neq \emptyset$ for $i=1, 2$, and $\alpha\cap S_{i}\neq\emptyset$ for $i=3, 4$. Hence $\psi_{F_{i}}(\gamma)\neq \emptyset$ for $i=1, 2$, and $\psi_{F_{i}}(\alpha)\neq \emptyset$ for $i=3, 4$; where $\psi$ is defined in Section 3. \vskip 2mm

{\bf Assumption 2.} (1) Let $\delta$ be an essential simple closed curve on the torus $F_{1}$ such that $d_{\mathcal {C}(F_{2})}(\psi_{F_{2}}(\gamma), \delta)\geq 5$

(2) Let $X$ be a full complex of $\mathcal{C}(F_{2})$ such that $d_{\mathcal {C}(F_{2})}(\psi_{F_{2}}(\gamma), \mathcal {D}(H_{X}))\geq 24$, where $H_{X}$ is the handlebody obtained by attaching 2-handles to $F_{2}$ along the vertices of $X$ then 3-handles to capping off the spherical boundary components. \vskip 2mm

Let $V_{F_{2}}=V\cup H_{X}$, and $V_{F_{1}, F_{2}}$ be the handlebody obtained by doing a surgery on $V_{F_{2}}$ along the slope $\delta$ on $F_{1}$. By Assumption 1,  $g(S_{3})=1$, $g(S_{4})\geq 2$, $V_{F_{1}, F_{2}}$ is not a $I$-bundle over a compact surface with $S_{i}$ as a horizontal boundary for $i=3, 4$. By Lemma 2.4, $diam_{\mathcal{C}(S_{i})}(\pi_{S_{i}}(\mathcal{D}(V_{F_{1}, F_{2}})))\leq 12$ for $i=3, 4$. \vskip 2mm

{\bf Assumption 3.} (1) Let $r$ be an essential simple closed curve on the torus $F_{3}$ such that $d_{\mathcal {C}(F_{3})}(\psi_{F_{3}}(\mathcal{D}(V_{F_{1}, F_{2}})), r)\geq 24$.

(2) Let $Y$ be a full complex of $\mathcal{C}(F_{4})$ such that $d_{\mathcal {C}(F_{4})}(\psi_{F_{4}}(\mathcal{D}(V_{F_{1}, F_{2}})), \mathcal {D}(H_{Y}))\geq 24$, where $H_{Y}$ is the handlebody obtained by attaching 2-handles to $F_{4}$ along the vertices of $Y$ then 3-handles to capping off the spherical boundary components. \vskip 2mm

Let $W_{F_{4}}=W\cup H_{Y}$, and $W_{F_{3}, F_{4}}$ be the handlebody obtained by doing a surgery on $W_{F_{4}}$ along the slope $r$ on $F_{3}$. Now both $M^{*}=V_{F_{2}}\cup_{S} W_{F_{4}}$ and $V_{F_{1}, F_{2}}\cup_{S} W_{F_{3}, F_{4}}$ are Heegaard splittings. Furthermore, we can prove that these two Heegaard splittings have distance 2  by  using Lemma 2.2 to take place of Lemma 2.3 in the proof of Proposition 3.1.

%Now let $X$ be a full complex of $\mathcal{C}_{2})$ such that $d_{\mathcal {C}(F_{2})}(\psi_{F_{2}}(\gamma), \mathcal {D}(H_{X}))\geq 24$, where $\mathcal {D}(H_{X})$ is the handlebody obtained by attaching handles to $F_{2}$ along the vertices of $X$. Similarly, let $Y$ be a full complex of $\mathcal{C}(F_{4})$ such that $d_{\mathcal {C}(F_{4})}(\psi_{F_{4}}(\alpha), \mathcal {D}(H_{Y}))> 24$.  Let $V_{F_{2}}$ be the compression body obtained by attaching the Handle body $H_{X}$ to $V$, and $W_{F_{4}}$ be the handlebody obtained by attaching $H_{Y}$ to $W$. Then $M^{*}=V_{F_{2}}\cup_{S} W_{F_{4}}$ is a Heegaard splitting. We can prove that the distance of $V_{F_{2}}\cup_{S} W_{F_{4}}$ is 2 Hence $M^{*}=V_{F_{2}}\cup_{S} W_{F_{4}}$ is irreducible and $\partial$-irreducible. \vskip 3mm

\vskip 2mm

%Then we can attach 2-handles along $\alpha$ and $\gamma$. Hence we get a Heegaard splitting
%$V\cup_{S} W$,  which has distance 2. And typically, $\partial _{-}V=F_{1}\cup F_{2}$ where $g(F_{1})=2$,
%$\partial _{-} W=F_{3}\cup F_{4}$ where $g(F_{3})=2$. Since $\gamma\cap S_{1}\neq \emptyset$,
%$\psi_{F_{1}}(\gamma)\neq \emptyset$. By theorem 2.7, there is
% a full simplex $X$ on $F_{1}$ such that $d_{\mathcal {C}(F_{1})}(\psi_{F_{1}}(\gamma), \mathcal {D}(H_{X}))\geq 24$. And we attach 2-handles
% along $X$. Then the compressionbody $V$ is changed into $V_{F_{1}}$. Hence it generates a new Heegaard splitting
% $V_{F_{1}}\cup_{S} W$. It is not hard to see that the Heegaard distance is still 2. And it is not hard to see that
 %any essential disk in $V_{F_{1}}$ can be isotopic to be disjoint from disk bounds by $C$.\vskip 3mm

 %On the other side, we can also find a full simplex $Y$ such that $d_{\mathcal {C}(F_{3})}(\mathcal {D}(H_{Y}), \psi_{F_{3}}(\alpha))\geq 24$.
 %After attaching 2-handles along $Y$, $W$ is changed into $W_{F_{3}}$. And the new Heegaard splitting
 %$V_{F_{1}}\cup_{S} W_{F_{3}}$ also has distance 2.\vskip 3mm

 %Since Heegaard distance of  $V_{F_{1}}\cup_{S} W_{F_{3}}$ is 2, $M=V_{F_{1}}\cup_{S} W_{F_{3}}$ is irreducible and $\partial$-irreducible. It means $M$ is a Haken 3-manifold. Now we start to prove that there is no essential torus or annulus in $M$.\vskip 3mm
Now we consider $M^{*}=V_{F_{2}}\cup_{S} W_{F_{4}}$. Note that $M^{*}$ has only two toral components. Since the distance of $V_{F_{2}}\cup_{S} W_{F_{4}}$ is 2, $M^{*}$ is irreducible and $\partial$-irreducible. \vskip 2mm

{\bf Claim 1.}  $M^{*}$ is atoroidal.\vskip3mm

{\bf Proof.} Suppose, otherwise, that $M^{*}$ contains an essential torus $T$. Since the distance of $V_{F_{2}}\cup_{S} W_{F_{4}}$ is 2,  $V_{F_{2}}\cup_{S} W_{F_{4}}$ is strongly irreducible. By Schultens's lemma, each component of $T\cap S$ is essential on both $T$ and $S$. Hence each component of $T\cap V_{F_{2}}$ and $T\cap W_{F_{4}}$ is an essential annulus in $V_{F_{2}}$ or $W_{F_{4}}$. \vskip 2mm

Let $A_{0}$ be one component of $T\cap V_{F_{2}}$. We first prove that there is one component of $\partial A_{0}$, say $a_{0}$, is not isotopic to $\beta$.

Now $V_{F_{2}}$ contains a $\partial$-compressing disk $B^{*}$ of $A_{0}$. By doing a surgery on $A_{0}$ along $B^{*}$, we can get a  disk $B_{0}$ in $V_{F_{2}}$. Since $A_{0}$ is essential, $B_{0}$ is essential. Suppose that the two components of $\partial A_{0}$ are isotopic to $\beta$. Since $\beta$ is non-separating on $S$, $\partial B_{0}$ bounds a once-punctured torus containing $\beta$, see Figure 8.

\begin{center}
\includegraphics[height=6cm, width=6cm]{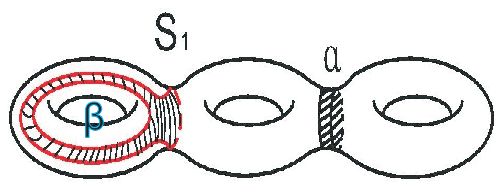}
\begin{center}
Figure 8
\end{center}
\end{center}
 By Assumption 1, $\beta\subset S_{2}$. Since $S_{2}$ has genus $g-1\geq 2$, $\partial B_{0}$ is not isotopic to $\alpha=\partial S_{2}$. By a standard outermost argument, $\psi_{F_{2}}(\partial B_{0})$ bounds an essential disk in $H_{X}$. Therefore $d_{\mathcal {C}(F_{2})}(\mathcal {D}(H_{X}), \psi_{F_{2}}(\beta))\leq 1$. Since $\gamma\cap \beta=\emptyset$, $d_{\mathcal {C}(F_{2})}(\psi_{F_{2}}(\beta), \psi_{F_{2}}(\gamma))\leq 1$. Hence $d_{\mathcal {C}(F_{2})}(\mathcal {D}(H_{X}), \psi_{F_{2}}(\gamma))\leq 2$. It contradicts Assumption 2. \vskip 3mm

Let $A_{1}$ be a component of $T\cap W_{F_{4}}$ which is incident to $A_{0}$. This means  that $a_{0}$ is one component of $\partial A_{1}$.
\vskip 3mm

{\bf Case 1.} \  $a_{0}\cap \alpha=\emptyset$, and $a_{0}\cap \gamma=\emptyset$.

Recall the definition of the surface $S^{\beta}$. Since $a_{0}$ is not isotopic to $\beta$,  $a_{0}\cap S^{\beta}\neq \emptyset$.  Since $\alpha, \gamma\subset S^{\beta}$, $d_{\mathcal {C}(S_{\beta})}(\pi_{S^{\beta}}(a_{0}), \alpha)\leq 1$, and $d_{\mathcal {C}(S^{\beta})}(\gamma, \pi_{S^{\beta}}(a_{0}))\leq 1$. Hence $d_{\mathcal {C}(S^{\beta})}(\alpha, \gamma)\leq 2$. This contradicts Assumption 1. \vskip 3mm

{\bf Case 2.}   $a_{0}\cap (\alpha\cup \gamma)\neq\emptyset$.

We assume that $a_{0}\cap \alpha\neq \emptyset$. By the above argument, $B_{0}$ is an essential disk in $V_{F_{2}}$ such that $\partial B_{0}$ is disjoint from $a_{0}$. Furthermore, $\partial B_{0}$ is not isotopic to $\alpha$. Since $B$ cuts $V_{F_{2}}$ into $F_{1}\times I$ and a handlebody $H$ such that $S_{2}\cup B=\partial H$, $\partial B_{0}\cap S_{2}\neq \emptyset$. Furthermore, all outermost disks of $B_{0}\cap B$ on $B_{0}$ lie in $H$. Hence $\pi_{S_{2}}(\partial B_{0})$ bounds an essential disk in $H$. This means $\psi_{F_{2}}(\partial B_{0})$ bounds an
essential disk in $H_{X}$.

If $a_{0}\cap \gamma=\emptyset$, then

\begin{center}
$d_{\mathcal {C}(F_{2})}(\psi_{F_{2}}(\partial B_{0}), \psi_{F_{2}}(\gamma))\leq d_{\mathcal {C}(F_{2})}(\psi_{F_{2}}(\partial B_{0}), \psi_{F_{2}}((a_{0}))+d_{\mathcal {C}(F_{2})}(\psi_{F_{2}}(a_{0}), \psi_{F_{2}}(\gamma))\leq 4$.
\end{center}

It contracts Assumption 2. Hence $a_{0}\cap \gamma\neq \emptyset$, and
$\psi_{F_{4}}(a_{0})\neq \emptyset$.\vskip 3mm

%{\bf Case 2.1.}  $a_{1}\cap \gamma= \emptyset$.

%We can assume that $a_{1}\cap \alpha \neq \emptyset$. For if $a_{1}\cap \alpha=\emptyset$, we can get that $d_{\mathcal {C}(S_{\beta})}(\alpha, \gamma)\geq 2$. It contradicts the choice of
%$\beta$. Then by the argument in previous chapter, we can get that $d_{\mathcal {C}(F_{1})}(\psi_{F_{1}}(b), \gamma)\leq 8$. It contracts the choice of $X$. End.
Since $A_{1}$ is an essential annulus in $W_{F_{4}}$, there is an essential disk $D_{0}$ obtained by doing boundary
compression on $A_{1}$ in $W_{F_{4}}$. Even more $\partial D_{0}\cap a_{0}=\emptyset$. Since $D$ cuts $W_{F_{4}}$ into  $F_{3}\times I$ and a handlebody $H^{*}$ containing $H_{Y}$, all outermost disks of $D_{0}\cap D$ in $D_{0}$ lies in $H^{*}$. Hence $\psi_{F_{4}}(\partial D_{0})$ bounds an essential disk in $H_{Y}$. Hence $\pi_{S_{4}}(\partial D_{0})\neq \emptyset$. Since $\partial D_{0}\cap a_{0}=\emptyset$, by Lemma 2.2, $d_{\mathcal {C}(S_{4})}(\pi_{S_{4}}(\partial D_{0}), \pi_{S_{4}}(a_{0}))\leq 2$. According to the definition of $\psi_{F_{4}}$, $d_{\mathcal {C}(F_{4})}(\psi_{F_{4}}(\partial D_{0}), \psi_{F_{4}}(a_{0}))\leq 2$.

Recall that the essential disk $B_{0}$ is obtained by doing a surgery on $A_{0}$ along a $\partial$-compressing disk in $V_{F_{2}}$. Since the distance of $V_{F_{2}}\cup_{S} W_{F_{4}}$ is two, $\partial B_{0}\cap \gamma\neq \emptyset$. Since $g(S_{3})=1$ and $g(S_{4})\geq 2$,
$V_{F_{2}}$ is not  a I-bundle of compact surface with $S_{4}$ as one horizontal boundary, by Lemma 2.4, $d_{\mathcal {C}(S_{4})}(\pi_{S_{4}}(\partial B_{0}), \pi_{S_{4}}(\alpha))\leq 12$. Hence $d_{\mathcal {C}(F_{4})}(\psi_{F_{4}}(\partial B_{0}), \psi_{F_{4}}(\alpha))\leq 12$. Since $\partial B_{0}\cap a_{0}=\emptyset$, $d_{\mathcal {C}(F_{4})}(\psi_{F_{4}}(\partial B_{0}), \psi_{F_{4}}(a_{0}))\leq 2$. It means that

\begin{center}

$d_{\mathcal {C}(F_{4})}(\psi_{F_{4}}(\partial D_{0}), \psi_{F_{4}}(\alpha))\leq d_{\mathcal {C}(F_{4})}(\psi_{F_{4}}(\partial D_{0}), \psi_{F_{4}}(a_{0}))+d_{\mathcal {C}(F_{4})}(\psi_{F_{4}}(\partial B_{0}), \psi_{F_{4}}(a_{0}))+d_{\mathcal {C}(F_{4})}(\psi_{F_{4}}(\partial B_{0}), \psi_{F_{4}}(\alpha))\leq 16$.
\end{center}     It contradicts Assumption 3.  END(Claim 1)\vskip 3mm

{\bf Claim 2.}  $M^{*}$ is anannular.\vskip 3mm

{\bf Proof.} Suppose, otherwise, that $M^{*}$ contains an essential torus $A$. Since the distance of $M^{*}=V_{F_{2}}\cup_{S} W_{F_{4}}$ is 2,  $M^{*}=V_{F_{2}}\cup_{S} W_{F_{4}}$ is strongly irreducible. By Schultens's lemma, each component of $A\cap S$ is essential on both $A$ and $S$. Hence each component of $A\cap V_{F_{2}}$ and $A\cap W_{F_{4}}$ is  either a spanning annulus or  an essential annulus with two boundary components lying on $S$.  There are four cases:\vskip 3mm

{\bf Case 1.} $|A\cap S|\geq 4$.

In this case, let $A_{0}$ be one component of $A\cap V_{F_{2}}$ and $A\cap W_{F_{4}}$ such that each of the two components of $A\cap V_{F_{2}}$ and $A\cap W_{F_{4}}$ incident to $A_{0}$ has its two boundary components lying on $S$. By the proof of Claim 1, the claim holds.\vskip 3mm

{\bf Case 2.} $|A\cap S|=1$.

Now $A$ intersects $S$ in an essential simple closed curve $a$. Furthermore, $a$, together with an essential simple closed curve $c_{1}$ on $F_{1}$, bounds a spanning annulus $A_{1}$ in $V_{F_{2}}$, and $a$, together with an essential simple closed curve $c_{2}$ on $F_{3}$, bounds a spanning annulus $A_{2}$ in $W_{F_{4}}$. Hence $a=c_{1}$ in $H_{1}(V_{F_{2}})$. Since $B$ cuts $V_{F_{2}}$ into $F_{1}\times I$ and a handlebody $H$ containing $\beta$, $a\neq \beta$ in $H_{1}(V_{F_{2}})$. Hence $a$ is not isotopic to $\beta$. There are three sub-cases: \vskip 3mm

{\bf Case 2.1.} $a\cap \alpha=\emptyset$, and $a\cap\gamma=\emptyset$.

In this case, $a\subset S_{1}, S_{3}$, and $\beta\subset S_{2}, S_{4}$, Hence $a\cap \beta=\emptyset$. This means that $d_{\mathcal{C}(S^{\beta})}(\alpha, \gamma)\leq 2$. It contradicts Assumption 1. \vskip 3mm

{\bf Case 2.2.} $a\cap \alpha=\emptyset$, and $a\cap\gamma\neq \emptyset$.

Now $A_{2}\cap D\neq\emptyset$. Let $c$ be an outermost arc of $A_{2}\cap D$ on $D$. This means that $c$, together with a sub-arc $c^{*}$ of $\gamma=\partial D$, bounds a disk $D^{*}$ in $D$ such that $D^{*}\cap D=c^{*}$. Now we can obtain a disk $D_{0}$ by doing a surgery on $A_{2}$ along $D^{*}$, say $D_{0}$. Hence $\partial D_{0}\cap a=\emptyset$. Furthermore, $D_{0}$ is an essential disk in $W_{F_{4}}$. Otherwise, we can reduce $|A_{2}\cap \gamma|$. Since $D$ cuts $W_{F_{4}}$ into $F_{3}\times I$ and a handlebody $H^{*}$ containing $S_{4}$. Hence $\partial D_{0}\cap S_{4}\neq\emptyset$. Furthermore, one outermost disk of $D_{0}\cap D$ on $D_{0}$ lies in $H^{*}$. Otherwise, we can reduce $|A_{2}\cap \gamma|$.  Hence $\psi_{F_{4}}(\partial D_{0})$ bounds an essential disk in $H_{Y}$.  By the assumption, $a\cap \gamma\neq \emptyset$. Hence $a\cap S_{4}\neq\emptyset$. Since the distance of $V_{F_{2}}\cup_{S} W_{F_{4}}$ is 2,
$\alpha\cap S_{4}\neq\emptyset$. Now by Lemma 2.2 and the above argument,
\begin{center}
$d_{\mathcal{C}(F_{4})}(\psi_{F_{4}}(\partial D_{0}), \psi_{F_{4}}(\alpha))\leq d_{\mathcal{C}(F_{4})}(\psi_{F_{4}}(\partial D_{0}), \psi_{F_{4}}(a))+d_{\mathcal{C}(F_{4})}(\psi_{F_{4}}(\alpha), \psi_{F_{4}}(a))\leq 4$.

\end{center}  It contradicts Assumption 3. \vskip 3mm

{\bf Case 2.3.} $a\cap \alpha\neq\emptyset$, and $a\cap\gamma\neq \emptyset$.

Let $D_{0}$ be as in Case 2.2. Since $a\cap \alpha\neq \emptyset$, $A_{1}\cap B\neq\emptyset$. Let $b$ be an outermost arc of $A_{1}\cap B$ on $B$. This means that $b$, together with a sub-arc $b^{*}$ of $\alpha$, bounds a disk $B^{*}$ in $B$ such that $B^{*}\cap B=b^{*}$. Now we can obtain a disk $B_{0}$ by doing a surgery on $A_{1}$ along $B^{*}$, say $B_{0}$. Hence $\partial B_{0}\cap a=\emptyset$. Similarly, $B_{0}$ is an essential disk in $V_{F_{2}}$. Since the distance of $V_{F_{2}}\cup_{S} W_{F_{4}}$ is 2, $B_{0}\cap S_{4}\neq\emptyset$.
By Lemma 2.4, $d_{\mathcal{C}(S_{4})}(\pi_{S_{4}}(\partial B_{0}), \pi_{S_{4}}(\alpha))\leq 12$. By Lemma 2.2,  $d_{\mathcal{C}(S_{4})}(\pi_{S_{4}}(\partial B_{0}), \pi_{S_{4}}(a))\leq 2$. By the argument in Case 2.2, $d_{\mathcal{C}(S_{4})}(\pi_{S_{4}}(\partial D_{0}), \pi_{S_{4}}(a))\leq 2$. Now  we have
\begin{center}
$d_{\mathcal{C}(F_{4})}(\psi_{F_{4}}(\partial D_{0}), \psi_{F_{4}}(\alpha))\leq d_{\mathcal{C}(S_{4})}(\pi_{S_{4}}(\partial D_{0}), \pi_{S_{4}}(\alpha))\leq 16.$
\end{center} Note that $\psi_{F_{4}}(\partial D_{0})$ is an essential disk in $H_{Y}$. It contradicts Assumption 3. \vskip 3mm

% Let $A_{0}$ be one component of $T\cap V_{F_{2}}$. We first prove that one component of $\partial A_{0}$, say $a_{0}$, is not isotopic to $\beta$.

%Now $V_{F_{2}}$ contains a $\partial$-compressing disk $B^{*}$ of $A_{0}$. By doing a surgery on $A_{0}$ along $B^{*}$, we can get a  disk $B_{0}$ in $V_{F_{2}}$. Since $A_{0}$ is essential, $B_{0}$ is essential. Otherwise, we can reduce $|T\cap S|$. Suppose that the two components of $\partial A_{0}$ are isotopic to $\beta$. Since $\beta$ is non-separating on $S$. Hence $\partial B_{0}$ bounds a once-punctured torus containing $\beta$.  By the assumptions, $\beta\subset S_{2}$, and $S_{2}$ has genus $g-1$, $\partial B_{0}$ is  not isotopic to $\alpha=\partial S_{2}$.   $\psi_{F_{2}}(\partial B_{0})$ bounds an essential disk in $H_{X}$. Since $\beta\subset S_{2}$, $d_{\mathcal {C}(F_{2})}(\mathcal {D}(H_{X}), \psi_{F_{2}}(\beta))\leq 1$. Since $\gamma\cap \beta=\emptyset$, $d_{\mathcal {C}(F_{2})}(\psi_{F_{2}}(\beta), \psi_{F_{2}}(\gamma))\leq 1$. Hence $d_{\mathcal {C}(F_{2})}(\mathcal {D}(H_{X}), \psi_{F_{2}}(\gamma))\leq 2$. This contradicts the choice of $X$. \vskip 3mm

{\bf Case 3.} $|A\cap S|=2$.

Now we may assume that $A\cap V_{F_{2}}$ contains two spanning annulus, and $A\cap W_{F_{2}}$ is an annulus with its two boundary components lying on $S$. By the arguments in Claim 1 and Case 2, Claim 2 holds. \vskip 3mm

{\bf Case 4.} $|A\cap S|=3$.

This case immediately from Claim 1 and Case 3. END(Claim 2) \vskip 3mm

Now $M^{*}$ is a hyperbolic 3-manifold, $M^{*}=V_{F_{2}}\cup_{S} W_{F_{4}}$ is a distance 2 Heegaard splitting of genus $g$. Furthermore, $M^{*}$ contains two toral boundary components $F_{1}$ and $F_{3}$. By the main results in [1] and  [16],  there are at most ten slopes $\delta$ on $F_{1}$ such that the manifold $M^{*}(\delta)$ obtained by doing Dehn filling on $M^{*}$ along $\delta$ is non-hyperbolic. By Assumption 2, there are infinitely many slopes $\delta$ so that $M^{*}(\delta)$ has a distance 2 Heegaard splitting of genus $g$. Hence there is at least one slope $\delta$ on $F_{1}$ such that  $M^{*}(\delta)$ is hyperbolic and $M^{*}(\delta)$ admits a distance 2 Heegaard splitting of genus $g$.  Similarly, by Assumption 3, there is a hyperbolic closed manifold which admits a distance 2 Heegaard splitting of genus $g$. END \vskip 3mm

\vskip 3mm

{\bf Acknowledgements} The authors thank Mario Eudave-Munoz and Jiming Ma for some helpful discussions.  \vskip 3mm

Ruifeng Qiu, Department of Mathematics, East China Normal University, Dongchuan Road 500, Shanghai 200241, China

rfqiu@math.ecnu.edu.cn

Yanqing Zou, School of Mathematical Sciences, Dalian University of Technology, Dalian 116022, China

yanqing\_dut@163.com

Qilong Guo,  School of Mathematical Sciences, Dalian University of Technology, Dalian 116022, China

guolong1999@yahoo.com.cn

\end{document}